\documentclass[opre,nonblindrev]{informs3} 
\usepackage{hyperref}
\DoubleSpacedXI 

\usepackage{endnotes}
\let\footnote=\endnote
\let\enotesize=\normalsize
\def\notesname{Endnotes}%
\def\makeenmark{$^{\theenmark}$}
\def\enoteformat{\rightskip0pt\leftskip0pt\parindent=1.75em
  \leavevmode\llap{\theenmark.\enskip}}

\usepackage{algpseudocode,multirow, float,bm,amsmath, theorem,accents}
\newcommand{\md}[1]{\mathbb{#1}}
\newcommand{\ma}[1]{\mathcal{#1}}

\newcommand{\thickhline}{%
	\noalign {\ifnum 0=`}\fi \hrule height 1pt
	\futurelet \reserved@a \@xhline
}

\renewcommand{\d}{\mathrm{d}}

\newtheorem{prop}{Proposition}
\newtheorem{theo}{Theorem}
\newtheorem{lem}{Lemma}
\newtheorem{defi}{Definition}
\newcommand*{\QED}{\hfill\ensuremath{\square}}
\usepackage[toc,page]{appendix}

\renewcommand{\Comment}[2][.55\linewidth]{%
	\leavevmode\hfill\makebox[#1][l]{\#~#2}}

\DeclareMathOperator\diag{Diag}
\DeclareMathOperator\rank{Rank}

\DeclareMathOperator\erf{erf}
\usepackage{algorithm}
\usepackage{natbib}
 \bibpunct[, ]{(}{)}{,}{a}{}{,}%
 \def\bibfont{\small}%
\TheoremsNumberedThrough     
\ECRepeatTheorems

\EquationsNumberedThrough    

\begin{document}


\RUNAUTHOR{Bertsimas and Li}

\RUNTITLE{Interpretable Matrix Completion}

\TITLE{Interpretable Matrix Completion: A Discrete Optimization Approach}

\ARTICLEAUTHORS{%
\AUTHOR{Dimitris Bertsimas}
\AFF{Sloan School of Management and Operations Research Center, 
	Massachusetts Institute of Technology, 
	Cambridge, MA 02139, \EMAIL{dbertsim@mit.edu}} 
\AUTHOR{Michael Lingzhi Li}
\AFF{Operations Research Center, 
	Massachusetts Institute of Technology, 
	Cambridge, MA 02139,  \EMAIL{mlli@mit.edu}}
} 

\ABSTRACT{%
We consider the problem of matrix completion on an $n \times m$ matrix. We introduce the problem of \emph{Interpretable Matrix Completion} that aims to provide meaningful insights for the low-rank matrix using side information. We show that the problem can be reformulated as  a binary convex optimization  problem. We design  OptComplete, based on a novel concept of stochastic cutting planes to enable efficient scaling of the algorithm up to matrices of sizes $n=10^6$
and $m=10^6$. We report 
experiments on both synthetic and real-world datasets that show that OptComplete has favorable scaling behavior and accuracy when compared with state-of-the-art methods for other types of matrix completion, while providing  insight on the factors that affect the matrix. 
}%


\KEYWORDS{Matrix Completion, Mixed-Integer Optimization,  Stochastic Approximation} 

\maketitle

%

\section{Introduction}
Low-rank matrix completion has attracted much attention after the successful application in the Netflix Competition. It is now widely utilized in far-reaching areas such as computer vision (\cite{candes2010matrix}), signal processing (\cite{ji2010robust}), and control theory (\cite{boyd1994linear}) to generate a completed matrix from partially observed entries. 

The classical low-rank matrix completion problem considers the following problem: Given a matrix $\bm{A}\in \md{R}^{n \times m}$ with entries only partially known (denote $\Omega \subset \{1,\ldots, n\} \times \{1,\ldots, m\}$ as the set of known entries), we aim to recover a matrix $\bm{X} \in \md{R}^{n \times m}$ of rank $k$ that minimizes a certain distance metric between $\bm{X}$ and $\bm{A}$ on the known entries of $\bm{A}$:
\begin{equation*}
\min_{\bm{X}} \frac{1}{nm} \sum_{(i,j) \in \Omega} \|X_{ij}-A_{ij}\|  \quad \text{subject to} \quad \rank(\bm{X})=k,
\end{equation*}
where we normalized the objective so that it is $O(1)$. The rank $k$ constraint on $\bm{X}$ can be equivalently formulated as the existence of two matrices $\bm{U}\in \md{R}^{n \times k}$, $\bm{V} \in \md{R}^{m \times k}$ such that $\bm{X}=\bm{U}\bm{V}^T$. Therefore, the problem can be restated as:
\begin{equation}
\min_{\bm{U}} \min_{\bm{V}} \frac{1}{nm} \sum_{(i,j) \in \Omega} \|X_{ij}-A_{ij}\|  \quad \text{subject to} \quad \bm{X}=\bm{U}\bm{V}^T. \label{eq:mcgeneral}
\end{equation}
In many applications for matrix completion, it is customary for each row of the data to represent an individual and each column a product or an item of interest, and $A_{ij}$ being the response data of individual $i$ on item $j$. Therefore, the matrices $\bm{U}$ and $\bm{V}$ are commonly interpreted as the ``user matrix" and ``product matrix" respectively. 

Let us denote each row of $\bm{U}$ as $\bm{u}_i$ and each column as $\bm{u}^i$ (similarly for $\bm{V}$). Then, $\bm{u}^i$ ($\bm{v}^i$) represents a ``latent feature" for users (products), and in total there are $k$ latent features for users (products). The goal of matrix completion is thus to discover such latent features of the users and the products, so that the dot product of such features on user $i$ and product $j$, $X_{ij}= \bm{u}_i \cdot \bm{v}_j$, is the intended response of user $i$ on product $j$. 

While this interpretation is  intuitive, it does not   offer   insight on what the latent features of users and products mean. Inductive Matrix Completion, first considered in \cite{dhillon2013new}, aims to rectify such problem by asserting that  each of the $k$ latent features is a linear combination of $p>k$ known features. We focus on the \emph{one-sided information} case, where only one of the user/product matrix is subject to such constraint. This case is more relevant as features related to users are often fragmented and increasingly constrained by data privacy regulations, while features about products are easy to obtain. We would also have a short discussion later on why the two-sided information case is not interesting under the context of this paper.

Without loss of generality, we would assume the information is on the product matrix and denote the known feature matrix with $p$ features as $\bm{B} \in \md{R}^{m \times p}$. As an example, if the items are movies, then $\bm{b}^j$ represents feature $j$ for a movie (actors, budget, running time, etc), and each product feature $\bm{v}^{j}$ needs to be linear combination of such features. Mathematically, this translates to the constraint:
\begin{equation*}
    \bm{V}=\bm{B}\bm{S}
\end{equation*}
where $\bm{S} \in \md{R}^{p \times k}$. Therefore, the inductive version of the problem in (\ref{eq:mcgeneral}) can be written as:
\begin{equation}
\min_{\bm{U}} \min_{\bm{S}} \frac{1}{nm} \sum_{(i,j) \in \Omega} \|X_{ij}-A_{ij}\|  \quad \text{subject to} \quad \bm{X}=\bm{U}\bm{S}^T\bm{B}^T. \label{eq:mcinductive}
\end{equation}
Although the inductive version of the problem adds more interpretability to the product latent features, they are still far from fully interpretable. For example, if we take the items to be movies, and features to be running time, budget, box office, and number of top 100 actors, then the generated features could look like:
\begin{align*}
    &5.6 \times \text{running time} - 0.00067 \times \text{budget} + 12 \times \text{\# of Top 100 actors},\\
    &0.25 \times \text{box office} -5 \times \text{\# of Top 100 actors}.
\end{align*}
These features, although a linear combination of interpretable features, are not very interpretable itself due to the involvement of multiple factors with different units. Therefore, it cannot significantly help decision makers to understand ``what is important" about the product. Furthermore, the appearance of the same factor in multiple features with different signs (as shown above) further complicates any attempt at understanding the result.

Therefore, we argue that instead of supposing the $k$ features are linear combinations of the $p$ known features, we should assume that the $k$ features are \emph{selected} from the $p$ known features. This formulation alleviates the two previous problems mentioned: it guarantees the latent features to be interpretable (as long as the original features are), and it prevents any duplicating features in the selected $k$ latent features. We denote this \emph{Interpretable Matrix Completion}. We use the term interpretable,  as opposed to inductive, to highlight that our approach, like sparse linear regression, gives actionable insights on what are the important features of matrix $\bm{A}$. We note that Interpretable Matrix Completion is considerably harder than inductive or the classical matrix completion problem as it is a discrete problem and selecting $k$ out of $p$ factors is exponential in complexity.

In this paper, we show that the Interpretable Matrix Completion problem can be written as a mixed integer convex optimization problem. Inspired by  \cite{SparseReg} for sparse linear regression, we reformulate the interpretable matrix completion problem as a binary convex optimization problem. Then we introduce a new algorithm OptComplete, based on stochastic cutting planes, to enable scalability for matrices of sizes on the order of $(n,m)=(10^6,10^6)$. In addition, we provide empirical evidence on both synthetic and real-world data that OptComplete is able to match or exceed current state-of-the-art methods for inductive and general matrix completion on both speed and accuracy, despite OptComplete solving a  more difficult problem.  

Specifically, our contributions in this paper are as follows:
\begin{enumerate}
	\item We introduce the interpretable matrix completion problem, and reformulate it as a binary convex optimization problem that can be solved using cutting planes methods. 
	\item We propose a new novel approach to cutting planes by introducing stochastic cutting planes. We prove that the new algorithm converges to an
	optimal solution of the interpretable matrix completion problem with exponentially vanishing failure probability. 
	\item We present computational results on both synthetic and real datasets that  show that the algorithm matches or outperforms current state-of-the-art methods in terms of both scalability and accuracy.
\end{enumerate}

The structure of the paper is as follows.
In Section  \ref{sec:MatrixCompletion}, we introduce the binary convex reformulation of the low-rank interpretable matrix completion problem, and how it can be solved through a cutting plane algorithm, which we denote CutPlanes.
In Section \ref{sec:OptComplete}, we introduce OptComplete, a stochastic  cutting planes method  designed to scale the CutPlanes algorithm in Section \ref{sec:MatrixCompletion}, and show that it recovers the optimal solution of CutPlanes with exponentially vanishing failure probability.
In Section \ref{sec:mc_compute_synt}, we report on computational  experiments with  synthetic data that compare OptComplete to  Inductive Matrix Completion (IMC) introduced in \cite{natarajan2014inductive} and SoftImpute-ALS (SIALS) by \cite{softimputeals}, two state-of-the-art matrix completion algorithms for inductive and general completion. We also compare OptComplete to CutPlanes to demonstrate the $20$x to $60$x speedup of the stochastic algorithm. In Section \ref{sec:mc_compute_real}, we report computational experiments on the Netflix Prize dataset. In Section \ref{sec_mc_conclusion} we provide our conclusions. 

\subsection*{Literature}
Matrix completion has been applied successfully to many tasks, including recommender systems \cite{koren2009matrix}, social network analysis \cite{chiang2014prediction} and clustering \cite{chen2014clustering}. After \cite{candesexact} proved a theoretical guarantee for the retrieval of the exact matrix under the nuclear norm convex relaxation, 
a lot of methods have focused on the nuclear norm problem (see \cite{SoftImpute}, \cite{beck2009fast}, \cite{jain2010guaranteed}, and \cite{tanner2013normalized} for examples).  Alternative methods include alternating projections by \cite{recht2013parallel} and Grassmann manifold optimization by \cite{keshavan2009matrix}. There has also been work where the uniform distributional assumptions required by the theoretical guarantees are violated, such as \cite{negahban2012restricted} and \cite{chen2014coherent}.

Interest in inductive matrix completion intensified after \cite{xuspeedup} showed that given predictive side information, one only needs $O(\log n)$ samples to retrieve the full matrix. 
Thus, most of this work (see \cite{xuspeedup}, \cite{jain2013provable}, \cite{farhat2013genomic},  \cite{natarajan2014inductive}) have focused on the case in which the side information is assumed to be perfectly predictive so that the theoretical bound of $O(\log n)$ sample complexity \cite{xuspeedup} can be achieved. 
\cite{chiang2015matrix} explored the case in which the side information is corrupted with noise, while \cite{shah2017matrix} and \cite{si2016goal} incorporated nonlinear combination of factors into the side information. Surprisingly, as pointed out by a recent article \cite{nazarov2018sparse}, there is a considerable lack of effort to introduce sparsity/interpretability into inductive matrix completion, with \cite{lu2016sparse}, \cite{soni2016noisy} and \cite{nazarov2018sparse} being among the only works that attempt to do so. Our work differs from the previous attempts in that previous attempts mainly focus on choosing latent features which are \emph{sparse} linear combinations of the given features. In contrast \emph{interpretable matrix completion} is aimed to \emph{select} exactly $k$ features from the known features. 
\section{Interpretable Matrix Completion}
\label{sec:MatrixCompletion}
In this section, we present the mathematical formulation of Interpretable Matrix Completion and how it can be reformulated as a binary convex problem that is based on \cite{SparseReg}. We show how this naturally leads to a cutting plane algorithm, and discuss its computational complexity. We also  discuss the two-sided information case, and how that reduces to the sparse regression problem.
\subsection{Binary Convex Reformulation of Interpretable Matrix Completion}
 The (one-sided) interpretable matrix completion problem can be written as a mixed binary optimization problem:
\begin{equation*}
\min_{\bm{U}} ~\min_{\bm{s}\in S_k^p} \frac{1}{nm} \sum_{(i,j) \in \Omega} \|X_{ij}-A_{ij}\|  \quad \text{subject to} \quad \bm{X}=\bm{U}\bm{S}\bm{B}^T,
\end{equation*}
where $\bm{S}=\diag\{s_1,\ldots, s_p\} \in \md{R}^{p \times p}$ and:
$$
S_{k}^p=\left\{\bm{s}=(s_1,\ldots, s_p)^T \in \{0,1\}^p:~ \sum_{i=1}^p s_i=k\right\}.
$$
We note that given that $\sum_{i=1}^p s_i=k$, the rank of matrix $\bm{X}$ is indeed $k$. We further note that the coefficients of $\bm{S}$ can be taken to be binary without loss of generality,
since if they are not and  $\bm{S}=\diag (1/d_1,\ldots, 1/d_p)$, then by applying 
the  transformation:
\begin{equation}
    \bm{U} \to \bm{U}\bm{D} \qquad \bm{S}\to \bm{S} \bm{D}^{-1} \label{eq:transform}
\end{equation}
for  $\bm{D}=\diag (d_1,\ldots, d_p)$, results in an equivalent problem with the coefficients of $\bm{S}$ being binary. 

For this paper, we consider the squared norm, and for robustness purposes (see \cite{SparseReg} and \cite{bc18}), we add a Tikhonov regularization term to the original problem. Specifically, the (one-sided) interpretable matrix completion problem with regularization we address is
\begin{equation}
\min_{\bm{U}} ~\min_{\bm{s}\in S_k^p}\frac{1}{nm}\left(\sum_{(i,j) \in \Omega} (X_{ij}-A_{ij})^2 + \frac{1}{\gamma}\|\bm{U}\|_2^2\right) \quad \text{subject to} \quad \bm{X}=\bm{U}\bm{S}\bm{B}^T \label{mainproblem}.
\end{equation}
In this section, we show how that problem (\ref{mainproblem}) can be reformulated as a binary convex optimization problem, and can be solved to optimality using a cutting plane algorithm. The main theorem and proof is presented below:
\begin{theo}
	\label{maintheo}
	Problem (\ref{mainproblem}) can be reformulated as a binary convex optimization problem:
	\begin{equation*}
	\min_{\bm{s}\in S_k^p} ~c(\bm{s})=\ \frac{1}{nm} \sum_{i=1}^n \overline{\bm{a}}_i \left(\bm{I}_m+\gamma \bm{W}_i\left(\sum_{j=1}^ps_j\bm{K}_j\right)\bm{W}_i\right)^{-1}\overline{\bm{a}}_i^T,
	\end{equation*}
	where $\bm{W}_1,\ldots, \bm{W}_n\in \md{R}^{m\times m}$
	are   diagonal  matrices:
	\[(\bm{W}_{i})_{jj}=\begin{cases}
	1, & (i,j) \in \Omega, \vspace{3pt} \\
	0, & (i,j) \not \in \Omega,
	\end{cases}\]
	$\overline{\bm{a}}_i=\bm{W}_i \bm{a}_i$, $i=1\ldots, n$, where $\bm{a}_i\in \md{R}^{1\times m}$ 
	is the $i$th row of $\bm{A}$ with unknown entries taken to be 0, 
	and $\bm{K}_j=\bm{b}^j(\bm{b}^j)^T \in \md{R}^{m\times m}$,  $j=1,\ldots, p$ with $\bm{b}^j \in \md{R}^{m\times 1}$
	the $j$th column of $\bm{B}$. 
\end{theo}
\proof{Proof:}
With the diagonal  matrices $\bm{W}_i$ defined above, we can rewrite the sum in (\ref{mainproblem}) over known entries of $\bm{A}$, $ \sum_{(i,j) \in \Omega} (X_{ij}-A_{ij})^2$, as a sum over the rows of $\bm{A}$:
\begin{equation*}
\sum_{i=1}^n \|(\bm{x}_i-\bm{a}_i)\bm{W}_i\|_2^2,
\end{equation*}
where $\bm{x}_i \in \md{R}^{1 \times m}$ is the $i$th row of $\bm{X}$. Using $\bm{X}=\bm{U}\bm{S}\bm{B}^T$, then  $\bm{x}_i=\bm{u}_i\bm{S}\bm{B}^T$ where  
$\bm{u}_i\in \md{R}^{1\times m}$ is the $i$th row of $U$. Moreover, 
$$\|\bm{U}\|_2^2=\sum_{i=1}^n \|\bm{u}_i\|_2^2.$$
Then, Problem (\ref{mainproblem}) 
becomes:
$$
\min_{\bm{s}\in S_k^p} ~\min_{\bm{U}} \frac{1}{nm}\left( \sum_{i=1}^n \left(\|
(\bm{u}_i\bm{S}\bm{B}^T-\bm{a}_i)\bm{W}_i\|_2^2 + \frac{1}{\gamma}\|\bm{u}_i\|_2^2\right) \right).
$$
We then notice that within the sum $\sum_{i=1}^n$  each row of $\bm{U}$ can be optimized separately, leading to:
\begin{equation}
\min_{\bm{s}\in S_k^p} ~\ \frac{1}{nm}\left( \sum_{i=1}^n \min_{\bm{u}_i}\left(\|
(\bm{u}_i\bm{S}\bm{B}^T-\bm{a}_i)\bm{W}_i\|_2^2 + \frac{1}{\gamma}\|\bm{u}_i\|_2^2\right) \right).\label{intermediate}
\end{equation}
The inner optimization problem $\displaystyle \min_{\bm{u}_i}\|
(\bm{u}_i\bm{S}\bm{B}^T-\bm{a}_i)\bm{W}_i\|_2^2 + \frac{1}{\gamma}\|\bm{u}_i\|_2^2$ can be solved in closed form given $\bm{S}$, as it is a weighted linear regression problem with Tiknorov regularization, see \cite{SparseReg}. The closed form solution is:
\begin{equation}
\min_{\bm{u}_i}\|
(\bm{u}_i\bm{S}\bm{B}^T-\bm{a}_i)\bm{W}_i\|_2^2 + \frac{1}{\gamma}\|\bm{u}_i\|_2^2= \overline{\bm{a}}_i(\bm{I}_m+\gamma \bm{W}_i\bm{B}\bm{S}\bm{B}^T\bm{W}_i)^{-1}\overline{\bm{a}}_i^T\label{eq:closed_form_sol}  .
\end{equation}
So Problem (\ref{intermediate}) can be simplified to:
$$
\min_{\bm{s}\in S_k^p} ~ \frac{1}{nm}\left( \sum_{i=1}^n \overline{\bm{a}}_i(\bm{I}_m+\gamma \bm{W}_i\bm{B}\bm{S}\bm{B}^T\bm{W}_i)^{-1}\overline{\bm{a}}_i^T\right) \notag .
$$
Finally, we notice that 
$$\bm{B}\bm{S}\bm{B}^T=\sum_{j=1}^p  s_j\bm{b}^j(\bm{b}^j)^T=\sum_{j=1}^p s_j\bm{K}_j,$$ 
and we obtain the required expression. Since $K_j$ are positive semi-definite, and the inverse of positive semi-definite matrices is a convex function, the entire function is convex in $\bm{s}$. 
\QED \endproof

With Theorem \ref{maintheo}, our original problem can now be restated as:
\begin{equation}
\label{eq:mc_central}
\min_{\bm{s}\in S_k^p} ~c(\bm{s})=\ \frac{1}{nm} \sum_{i=1}^n \overline{\bm{a}}_i \left(\bm{I}_m+\gamma \bm{W}_i\left(\sum_{j=1}^ps_j\bm{K}_j\right)\bm{W}_i\right)^{-1}\overline{\bm{a}}_i^T.
\end{equation}
This can be solved utilizing the cutting plane algorithm first introduced by \cite{CuttingPlane}, summarized as Algorithm \ref{alg:mc_alg0}. 

\begin{algorithm}[ht] 
	\begin{algorithmic}[1]
		\Procedure{CUTPLANES}{$\bm{A},\bm{B}$}\Comment{masked matrix $\bm{A}$, and feature matrix $\bm{B}$}
		\State $t \gets 1$
		\State $\bm{s}_1 \gets \text{warm start}$ \Comment{Heuristic Warm Start}
		\State $\eta \gets 0 $ \Comment{Initialize feasible solution variable}
		\While{$\eta_t<c(\bm{s}_t)  $}\Comment{While the current solution is not optimal}
		\State $\bm{s}_{t+1}, \eta_{t+1} \gets \displaystyle \argmin_{\bm{s}\in S_k^p, \eta>0} \quad \eta \quad \text{s.t.} \quad \eta \geq c(\bm{s}_i)+\nabla c(\bm{s}_i)^T(\bm{s}-\bm{s}_i) \quad \forall i \in [t]$
		\State $t \gets t+1$
		\EndWhile
		\State $\bm{s}\gets \bm{s}_t$
		\State $i \gets 1$
		\For{$i<n$} \Comment{Fill each row $\bm{x}_i$ of final output matrix $\bm{X}$}
		\State $\bm{x}_i \gets \bm{B}^{\bm{s}}((\bm{B}^{\bm{s}})^T\bm{W}_i\bm{B}^{\bm{s}})^{-1}(\bm{B}^{\bm{s}})^T\overline{\bm{a}}_i^T$ \Comment{$\bm{B}^{\bm{s}}$ is submatrix of $\bm{B}$ with $\bm{s}$ columns}
		\EndFor
		\State \textbf{return} $\bm{X}$\Comment{Return the filled matrix $\bm{X}$}
		\EndProcedure
	\end{algorithmic}
	\caption{Cutting-plane algorithm for   matrix completion with side information. }
	\label{alg:mc_alg0}
\end{algorithm}

The cutting plane algorithm, at iteration $t$, adds a linear approximation of $c(\bm{s})$ at the current feasible solution $\bm{s}_t$ to the set of constraints:
\begin{equation}
\eta\geq c(\bm{s}_t)+\nabla c(\bm{s}_t)^T (\bm{s}-\bm{s}_t) \label{cuttingplaneform},
\end{equation}
and we solve the mixed-integer linear programming problem:
\begin{align*}
&\min_{\bm{s}\in S_k^p,\eta \geq 0} \eta \\
&\eta \geq  c(\bm{s}_t)+\nabla c(\bm{s}_t)^T(\bm{s}-\bm{s}_t), \quad i \in [t]
\end{align*}
to obtain $\bm{s}_{t+1}, \eta_{t+1}$. We see that $\eta_{t+1}$ is exactly the minimum value of the current approximation of $c(\bm{s})$, $c_t(\bm{s})$, defined below:
\begin{equation*}
    \eta_{t+1} = \min_{\bm{s}} \max_{i \in [t]} c(\bm{s}_t)+\nabla x(\bm{s}_t)^T(\bm{s}-\bm{s}_t) = \min_{\bm{s}} c_t(\bm{s}).
\end{equation*}
Since $c(\bm{s})$ is convex, the piecewise linear approximation $c_t(\bm{s})$ is an outer approximation ($c_t(\bm{s})<=c(\bm{s}) \; \forall \bm{s}$), so $\eta_t\leq c(\bm{s}_t)\; \forall t$. As the algorithm progresses, the set of linear approximations form an increasingly better approximation of $c(\bm{s})$, and $\eta_t$ increases with $t$. The algorithm terminates once $\eta_t$ does not further increase, as it implies the linear approximation shares the same minimum value as the true function $c(\bm{s})$, which is the desired value. 

Once the optimal solution $\bm{s}^*$ is reached, we can obtain the optimal $\bm{U}$ using the closed form solution in (\ref{eq:closed_form_sol}) and recover $\bm{X}$. In the next section, we discuss how this algorithm can be implemented in the context of $c(\bm{s})$ in (\ref{eq:mc_central}) and derive its computational complexity.

\subsection{Implementation and Computational Complexity of CutPlanes}
The computational complexity of the cutting plane comes from calculating $c(\bm{s})$ and its derivative $\nabla c(\bm{s})$. We first introduce the notations $\alpha_i(\bm{s}) \in \md{R}$ and $\bm{\gamma}_i(\bm{s}) \in \md{R}^{m \times 1}$.
\begin{equation}
\label{eq:mv_alpha}
\alpha_i(\bm{s})= \frac{1}{m}\overline{\bm{a}}_i \bm{\gamma}_i(\bm{s})=\frac{1}{m}\overline{\bm{a}}_i\left[\left(\bm{I}_m+\gamma \bm{W}_i\left(\sum_{j=1}^p s_j\bm{K}_j\right)\bm{W}_i\right)^{-1}\overline{\bm{a}}_i^T\right],~~i=1,\ldots,n.
\end{equation}
Then, the function $c(\bm{s})$ in  (\ref{eq:mc_central}) can be expressed as
\begin{equation}
c(\bm{s})=\frac{1}{n}\sum_{i =1}^n \alpha_i(\bm{s})=\frac{1}{n} \sum_{i=1}^n  \frac{\overline{\bm{a}}_i \bm{\gamma}_i(\bm{s})}{m}.\label{cfbreakdown} 
\end{equation}
To calculate the derivative $\nabla c(\bm{s})$, it is easier to utilize the expression in Theorem \ref{maintheo} and then utilize the chain rule. After some manipulations, we obtain
\begin{equation}
\nabla c(\bm{s})=\frac{1}{n}\sum_{i =1}^n -\frac{\gamma \left(\bm{B}^T\bm{W}_i \bm{\gamma}_i(\bm{s})\right)^2}{m}.\label{cfderivbreakdown} 
\end{equation}
Therefore, we would focus on calculating $\bm{\gamma}_i(\bm{s})$. First, by the Matrix Inversion Lemma (\cite{matrixinv}) we have
\begin{align}
    \bm{\gamma}_i(\bm{s})&=\left(\bm{I}_m+\gamma \bm{W}_i\left(\sum_{j=1}^p s_j\bm{K}_j\right)\bm{W}_i\right)^{-1}\overline{\bm{a}}_i^T\nonumber \\&=\left(\bm{I}_m-\bm{V}\left(\frac{\bm{I}_k}{\gamma}+\bm{V}^T\bm{W}_i\bm{V}\right)^{-1}\bm{V}^T\right)\overline{\bm{a}}_i^T\nonumber \\&=\left(\overline{\bm{a}}_i^T
-\bm{V}\left(\frac{\bm{I}_k}{\gamma}+\bm{V}^T\bm{W}_i\bm{V}\right)^{-1}\bm{V}^T\overline{\bm{a}}_i^T\right) \label{eq:mv_gamma1},
\end{align}
where $\bm{V}\in \md{R}^{m\times k}$ is the feature matrix formed by the $k$ columns of $\bm{B}$ such that $s_j=1$, and we have suppressed the dependency of $\bm{V}$ on $\bm{s}$ for notation ease. 
Note that in  order to compute  $\gamma_i(\bm{s})$ using Eq.  (\ref{eq:mv_alpha}) we need to invert an $m\times m$ matrix, while
from Eq. (\ref{eq:mv_gamma1})  we need to invert  a $k\times k$ matrix $\frac{\bm{I}_k}{\gamma}+\bm{V}^T\bm{W}_i\bm{V}$, which only requires $O(k^3)$ calculations. Furthermore, note that calculating $\overline{\bm{a}}_i\overline{\bm{a}}_i^T$ only requires $|\Omega_i|$ multiplications where $\Omega_i$ is the number of known entries in row $i$ of $\bm{A}$ as we do not need to multiply on the unknown entries. Similarly, we can compute $\bm{V}^T\overline{\bm{a}}_i^T$ in $|\Omega_i|k$ multiplications, and $\bm{V}^T\bm{W}_i\bm{V}$ in  $|\Omega_i|k^2$ multiplications.

Therefore, we can compute $\gamma_i(\bm{s})$ in floating point complexity of $O(|\Omega_i |k^2+k^3)$. Then to calculate $\overline{\bm{a}}_i \bm{\gamma}_i(\bm{s})$ in  (\ref{cfbreakdown}) and  $-\gamma \left(\bm{B}^T\bm{W}_i \bm{\gamma}_i(\bm{s})\right)^2$ (\ref{cfderivbreakdown}) only requires $O(|\Omega_i|)$ and $O(|\Omega_i|p)$ calculations respectively. Thus, the total complexity of generating a full cutting plane is:
\begin{equation}
\sum_{i=1}^n O(|\Omega_i|p+|\Omega_i|k^2+k^3) =  O(|\Omega|(p+k^2)+nk^3) .\label{eq:fullcomplexity}
\end{equation}

\subsection{Two-sided Information Case}
In this section, we briefly discuss the matrix completion problem under the two-sided information case, and how it reduces to the problem
of sparse linear regression. 
The two sided interpretable matrix completion problem with Tikhonov regularization can be stated as follows:
\begin{equation}
\min_{\bm{L}}  \frac{1}{nm}\left(\sum_{(i,j) \in \Omega} (X_{ij}-A_{ij})^2 +\frac{1}{\gamma} \|\bm{L}\|_2^2 \right) \quad \text{subject to} \quad \bm{X}=\bm{U}\bm{L}\bm{B}^T \quad \|\bm{L}\|_0=k \label{twosidedproblem},
\end{equation}
where $\bm{U} \in \md{R}^{n \times p_1}$ is a known matrix of $p_1$ features of each row,  $\bm{B} \in \md{R}^{m \times p_2}$ is a known matrix of $p_2$ features of each column, and $\bm{L} \in \md{R}^{p_1 \times p_2}$ is a sparse matrix that has $k$ nonzero  entries, ensuring that $\rank(\bm{X})\leq k$. 
We note that in Eq.   \eqref{twosidedproblem} we restrict the support of matrix $\bm{L}$ to be $k$, rather than forcing the entries of $\bm{L}$ to be binary. This is because unlike in the one-sided case, both $\bm{U}$ and $\bm{B}$ are known, so we cannot apply the scaling transformation in (\ref{eq:transform}).

We denote by $\bm{u}^i\in \md{R}^{n \times 1} $  the $i$th  column of $\bm{U}$ and $\bm{b}^j\in \md{R}^{m \times 1}  $   the $j$th column of $\bm{B}$. 
We introduce the  matrices $\bm{W}_i$
as in Theorem \ref{maintheo}. Using $\bm{X}=\bm{U}\bm{L}\bm{B}$, we can write 
$$X_{ij}=\sum_{q=1}^{p_1} 
\sum_{\ell=1}^{p_2} 
L_{q,\ell} D^{q,\ell}_{ij},$$ where $D^{q,\ell}_{ij}=(\bm{u}^q(\bm{b}^{\ell})^T)_{ij}$ is the $(i,j)$th entry of the matrix formed by multiplying $q$th column of $\bm{U}$ with $\ell$th column of $\bm{B}$. Then, Problem (\ref{twosidedproblem}) becomes:
\begin{equation}
\label{eq:mc_inter2}
\min_{\bm{L}}  \frac{1}{nm}\left(\sum_{(i,j) \in \Omega} \left ( \sum_{q=1}^{p_1} 
\sum_{\ell=1}^{p_2} 
L_{q,\ell}D^{q,\ell}_{ij}-A_{ij}\right)^2 +\frac{1}{\gamma} \|\bm{L}\|_2^2 \right) \quad \text{subject to}  \quad \|\bm{L}\|_0=k.
\end{equation}
As every $\bm{D}$ matrix is known, this becomes a sparse regression problem where there are $p_1p_2$ features to choose from (the $\bm{D}$ matrices), there are $|\Omega|$ samples (the $\bm{A}$ matrix), the sparsity requirement is $k$, the regression coefficients are $\bm{L}$, and we have Tikhonov regularization. Vectorizing $\bm{D}$, $\bm{L}$, and $\bm{A}$ reduces the problem back to the familiar form of sparse linear  regression, that can be solved by  the algorithm developed in \cite{SparseReg} at scale.

\section{OptComplete: The Stochastic Cutting Plane Speedup}
\label{sec:OptComplete}
In this section, we introduce OptComplete, a stochastic version of the cutting plane algorithm introduced in Section \ref{sec:MatrixCompletion}. We present theoretical results to show that the stochastic algorithm recovers the true optimal solution of the original algorithm with high probability without distributional assumptions. We also include a discussion on the dependence of such probability with various factors and its favorable theoretical computational complexity.
\subsection{Introduction of OptComplete}
In the previous section, we showed that through careful  evaluation, we can calculate a full cutting plane in $O(|\Omega|(p+k^2)+nk^3)$ calculations. However, in very high dimensions where $|\Omega|, n, m$ are extremely large, the cost of generating the full cutting plane is still prohibitive. Thus, we consider generating approximations of the cutting plane that would enable the algorithm to scale for  high values for  $n$ and $m$. Specifically, consider the cutting plane function in (\ref{cfbreakdown}), reproduced below:
\begin{equation*}
c(\bm{s})=\frac{1}{n} \sum_{i=1}^n \alpha_i(\bm{s}),
\end{equation*}
where:
\begin{equation*}
    \alpha_i(\bm{s})=\frac{1}{m}\left(\overline{\bm{a}}_i\overline{\bm{a}}_i^T
-\overline{\bm{a}}_i\bm{V}\left(\frac{\bm{I}_k}{\gamma}+\bm{V}^T\bm{W}_i\bm{V}\right)^{-1}\bm{V}^T\overline{\bm{a}}_i^T\right).
\end{equation*}
We approximate the inner term $\alpha_i(\bm{s})$ by choosing $1\leq f<m$ samples from $\{1,\ldots,m\}$ \emph{without replacement}, with the set denoted $F$. Then we formulate the submatrix $\bm{V}_F$ with such selected rows, and similarly with $\overline{\bm{a}}_{Fi}$. Then we calculate the approximation:
\begin{equation*}
    \alpha_i(\bm{s})\approx \alpha^F_i(\bm{s})=\frac{1}{f}\left(\overline{\bm{a}}_{Fi}\overline{\bm{a}}_{Fi}^T
-\overline{\bm{a}}_{Fi}\bm{V}_F\left(\frac{\bm{I}_k}{\gamma}+\bm{V}^T_F\bm{W}_i\bm{V}_F\right)^{-1}\bm{V}^T_F\overline{\bm{a}}_{Fi}^T\right).
\end{equation*}
Then we choose $1 \leq g < n$ samples from $\{1,\ldots,n\}$ without replacement, with the set denoted $G$. We can then calculate an approximation of $c(\bm{s})$ using the approximated $\alpha^F_i(\bm{s})$:
\begin{equation*}
    c(\bm{s}) \approx \tilde{c}_G^F(\bm{s})=\frac{1}{r} \sum_{i \in G}  \alpha^F_i(\bm{s}),
\end{equation*}
where the set $F$ is chosen independently for every row $i \in G$. Then the derivative of $\tilde{c}_G^F(\bm{s})$ is:
\begin{equation*}
    \nabla \tilde{c}_G^F(\bm{s})=\frac{1}{r} \sum_{i \in G}  \frac{(\bm{B}^T\bm{W}_i\bm{\gamma}_i^F(\bm{s}))^2}{f}.
\end{equation*}
Using such approximations, we can derive a stochastic cutting-plane algorithm, which we call OptComplete  presented as  Algorithm \ref{alg:mc_alg1}.

\begin{algorithm}[h] 
	\begin{algorithmic}[1]
		\Procedure{OptComplete}{$\bm{A},\bm{B}$}\Comment{masked matrix $\bm{A}$, and feature matrix $\bm{B}$}
		\State $t \gets 1$
		\State $\bm{s}_1 \gets \text{random initialization}$
		\State $\eta \gets 0 $ \Comment{Initialize feasible solution variable}
		\While{$\eta_t<c(\bm{s}_t)  $}\Comment{While the current solution is not optimal}
		\State $G \gets \text{$|g|$-sized random sample of $\{1,\ldots, n\}$ with replacement}$ 
		\For{$i \in G$} \Comment{Generate $F$ for each row in random sample}
		\State $F_i \gets \text{$|f|$-sized random sample of $\{1,\ldots, m\}$ with replacement}$ 
		\EndFor
		\State $\bm{s}_{t+1}, \eta_{t+1} \gets \displaystyle \argmin_{\bm{s}\in S_k^p, \eta>0} \quad \eta \quad \text{s.t.} \quad \eta \geq \tilde{c}^F_G(\bm{s}_i)+\nabla \tilde{c}^F_G(\bm{s}_i)^T(\bm{s}-\bm{s}_i) \quad \forall i \in [t]$
		\State $t \gets t+1$
		\EndWhile
		\State $\bm{s}\gets \bm{s}_t$
		\State $i \gets 1$
		\For{$i<n$} \Comment{Fill each row $\bm{x}_i$ of final output matrix $\bm{X}$}
		\State $\bm{x}_i \gets \bm{B}^{\bm{s}}((\bm{B}^{\bm{s}})^T\bm{W}_i\bm{B}^{\bm{s}})^{-1}(\bm{B}^{\bm{s}})^T\overline{\bm{a}}_i^T$ \Comment{$\bm{B}^{\bm{s}}$ is submatrix of $\bm{B}$ with $\bm{s}$ columns}
		\EndFor
		\State \textbf{return} $\bm{X}$\Comment{Return the filled matrix $\bm{X}$}
		\EndProcedure
	\end{algorithmic}
	\caption{Stochastic Cutting-plane algorithm for   matrix completion with side information. }
	\label{alg:mc_alg1}
\end{algorithm}

For this algorithm to work, we need the approximation $\tilde{c}^F_G(\bm{s})$ and its derivative to be close to the nominal values. Furthermore, the approximated cutting planes should not cutoff the true solution. In the next section, we show that OptComplete enjoys such properties with high probability. In Section \ref{subsec:sampling}, we discuss how to select the size of $f$ and $g$.
\subsection{Main Theoretical Results}
We would first show that the inner approximation is close to the true term with high probability:
\begin{theo}
    \label{theo:inner_theo}
    Let $\bm{A}$ be a partially known matrix, $\bm{B}$ a known feature matrix, and $\bm{W}_i$ as defined in Theorem \ref{maintheo}. Let $F$ be a random sample of size $f$ from the set $\{1,\ldots,m\}$, chosen without replacement. With probability at least $1-\epsilon$, we have
    
    \begin{align*}
        |\alpha_i(\bm{s})-\alpha_i^F(\bm{s})| &\leq \sqrt{\frac{Mk\log(\frac{k}{\epsilon})}{f}}, \qquad \forall i \in \{1,\ldots,m\} ,\quad \forall \bm{s} \in S_k^p, \vspace{3pt}\\
                \left\|\frac{(\bm{B}^T\bm{W}_i\bm{\gamma}_i(\bm{s}))^2}{m}-\frac{(\bm{B}^T\bm{W}_i\bm{\gamma}_i^F(\bm{s}))^2}{f}\right\|_2 &\leq \sqrt{\frac{M'(p+k)\log(\frac{k}{\epsilon})}{f}}, \qquad \forall i \in \{1,\ldots,m\}, \quad \forall \bm{s} \in S_k^p,
    \end{align*}
    where $M,M'$ are absolute constants.
\end{theo}
We see that, without assumptions on the data, the inner approximation for both the value and the derivative follows a bound with $O\left(\sqrt{\frac{(p+k)}{f}}\right)$ terms with very high probability. Furthermore, inverting the statements give that, for all $i \in \{1,\ldots,m\}$ and all $\bm{s} \in S_k^p$:
\begin{align*}
    \md{P}\left(|\alpha_i(\bm{s})-\alpha_i^F(\bm{s})| \geq \delta\right)&\leq k\exp\left(-\frac{f\delta^2}{Mk}\right) , \vspace{3pt} \\
           \md{P}\left( \left\|\frac{(\bm{B}^T\bm{W}_i\bm{\gamma}_i(\bm{s}))^2}{m}-\frac{(\bm{B}^T\bm{W}_i\bm{\gamma}_i^F(\bm{s}))^2}{f}\right\|_2 \geq \delta\right)&\leq k\exp\left(-\frac{f\delta^2}{M'(p+k)}\right) .
\end{align*}
So the failure probability drops off exponentially with increasing bound $\delta$, reflecting a Gaussian tail structure for $\alpha_i(\bm{s})$ and $\bm{B}^T\bm{W}_i\bm{\gamma}_i^F(\bm{s})$.  The proof is contained in Appendix \ref{app:inner_theo_proof}. 

Using this result, we are able to prove a tight deviation bound for the approximated cost function $c^F_G(\bm{s})$ and $\nabla c^F_G(\bm{s})$:
\begin{theo}
    \label{theo:outer_theo}
    Let $\bm{A}$ be a partially known matrix, $\bm{B}$ a known feature matrix, and $\bm{W}_i$ as defined in Theorem \ref{maintheo}. Let $G$ be a random sample of size $g$ from $\{1,\ldots,n\}$ chosen without replacement. Then for each $i \in G$, we let $F_i$ be a random sample of size $f$ from the set $\{1,\ldots,m\}$, all chosen without replacement. We have, with probability at least $1-\epsilon$:
    \begin{align*}
        |\tilde{c}^F_G(\bm{s})-c(\bm{s})|&\leq \sqrt{\frac{Ak\log\left(\frac{k}{\epsilon}\right)}{g}}, \qquad \forall \bm{s} \in S_k^p , \vspace{3pt} \\
                \|\nabla \tilde{c}^F_G(\bm{s})-\nabla c(\bm{s})\|_2&\leq \sqrt{\frac{B(p+k)\log\left(\frac{k}{\epsilon}\right)}{g}}, \qquad \forall \bm{s} \in S_k^p,
    \end{align*}
    where $A,B$ are absolute constants.
\end{theo}
Similar to the inner approximations,  $c^F_G(\bm{s})$ and $\nabla c^F_G(\bm{s})$ has Gaussian tails. Furthermore, the scaling here only depends on $g$ and not $f$: This shows that the error of the inner approximation is dominated by the outer sampling of the rows $G$. The proof is contained in Appendix \ref{app:outer_theo_proof}.

Then, using this result, we are able to prove our main result for OptComplete. We would first introduce a new definition:
\begin{defi}
	\label{convexityparamdefi}
The \textbf{convexity parameter}  $a$ of the cost function $c(\bm{s})$ is defined as the largest positive number for which the the following statement is true:
\begin{equation}
c(\bm{s})\geq c(\bm{s}_0)+\nabla c(\bm{s}_0)^T(\bm{s}-\bm{s}_0)+\frac{a^2}{2}(\bm{s}-\bm{s}_0)^T(\bm{s}-\bm{s}_0) \quad \forall\bm{s},
\bm{s}_0 \in S^p_k\;\; \forall  i \label{convexityparam}
\end{equation}
\end{defi}
We have the following proposition which shows that unless the cost function is degenerate (i.e. different sets of $k$ features doesn't change the solution), we always have a positive convexity parameter: 
\begin{prop}
	\label{prop:strong_convexity_cond}
	Assume that there does not exist $\bm{s}_1,\bm{s}_0 \in S^p_k$ such that $c(\bm{s}_1)=c(\bm{s}_0)$. Then $a>0$. 
\end{prop}
The proof is contained in Appendix \ref{app:strong_convexity_cond_proof}.  Now, we state our main theorem for OptComplete. 
\begin{theo}
	\label{theo:mc_proof}
	For the matrix completion problem (\ref{mainproblem}), let $\bm{B} \in \md{R}^{m \times p}$ be a known feature matrix, $\bm{A} \in \md{R}^{n \times m}$ a matrix with entries partially known, and OptComplete as defined in Algorithm \ref{alg:mc_alg1}. Assume that Problem  (\ref{mainproblem}) is feasible. Then, OptComplete terminates in a finite number of steps $C$, and finds an optimal solution of (\ref{mainproblem}) with probability at least $1-kC\exp\left(-\frac{Da^4g}{(p+k)}\right)$ where $D$ is an absolute constant independent of $C, f,g,k,p$, and $a$ is the convexity parameter of the functions $\tilde{\alpha}_i^s(\bm{s})$.
\end{theo}
The proof is contained in Appendix \ref{app:optcomplete_theo_proof}. This theorem shows that as long as the original problem is feasible, OptComplete is able to find the optimal solution of the original binary convex problem with exponentially vanishing failure probability that scales as $O\left(\exp\left(\frac{-g}{(p+k)}\right)\right)$. The theorem requires no assumptions on the data, and thus applies generally. We again note that the bound does not depend on $f$ and only on $g$: we would discuss how this would inform our selection of the size of $f$ and $g$ in the next section.
\subsection{Sampling Size and Computational Complexity}
\label{subsec:sampling}
To select an appropriate $f$ and $g$, we first note that \cite{candesexact} showed that to complete a square $N\times N$ matrix of 
 rank $k$, we need at least $O(kN\log N)$ elements. Assume an average known rate of $\alpha=\frac{|\Omega|}{mn}$ in the original matrix $\bm{A}$, the expected number of known elements under a sampling of $f$ and $g$ is $\alpha fg$. Using $N^2=mn$, we need that:
 \begin{equation}
     \label{eq:sampling_lower_bound}
     \alpha fg\geq c\cdot k\sqrt{nm}\log  (\sqrt{nm})
 \end{equation}
 for some constant $c$. Theorem \ref{theo:mc_proof} showed that the bound on failure probability scales with $O\left(\exp\left(\frac{-g}{(p+k)}\right)\right)$, and thus we cannot have $g$ too small. Using (\ref{eq:fullcomplexity}), the complexity of the cutting plane with $f$ and $g$ samples are:
 \begin{equation}
  O(\alpha fg(p+k^2)+gk^3).
  \label{eq:reducedcomplexity}
 \end{equation}
 Therefore, if we fix the expected known elements ($\alpha fg$) constant, it is more advantageous to select a smaller $g$, as $g$ scales with $k^3$. Thus, we set:
 \begin{equation}
f=\min\left(\frac{ck\sqrt{mn}\log (\sqrt{mn})}{\alpha\min (g_0, n)},m\right) , \qquad g=\min(g_0, n) , \label{eq:choosecutsize}
\end{equation}
Experimentally, we found $g_0=100$, $c=1$ to generate good results (the results were similar for $\frac{1}{2}\leq c \leq 2$). Therefore, by (\ref{eq:reducedcomplexity}), the approximated cutting plane has a computational complexity of:
\begin{equation}
    O\left(k\sqrt{mn}\log(\sqrt{mn})(p+k^2)\right).
\end{equation}
This scales in a square root fashion in $n$ and $m$, rather than linearly in $n$ and $m$ for the full cutting plane. This allows OptComplete to enjoy a considerable speedup compared to CutPlanes, as  
 demonstrated in Section \ref{sec:mc_compute_synt}. 
\section{Synthetic Data Experiments}
\label{sec:mc_compute_synt}
We assume that the matrix $\bm{A}=\bm{UV}+\bm{E}$, where $\bm{U} \in \md{R}^{n \times k}$, $\bm{V} \in \md{R}^{k\times m}$, and $\bm{E}$ is an error matrix with individual elements sampled from $N(0,0.01)$. We sample the elements of $\bm{U}$ and $\bm{V}$ from a uniform distribution of $[0,1]$, and then randomly select a fraction  $\mu=1-\alpha$ to be missing.
We formulate the feature matrix $\bm{B}$ by combining $\bm{V} \in \md{R}^{k\times m}$ with a confounding matrix $\bm{Z} \in \md{R}^{ (p-k)\times m }$ that contains unnecessary factors sampled similarly from the Uniform $[0,1]$ distribution.
We run   OptComplete on a server with $16$ CPU cores, using Gurobi 8.1.0. For each combination $(m,n,p,k,\mu)$, we ran 10 tests and report  the median value for every statistic.

We report the following statistics with $\bm{s}^*$ being the ground-truth factor vector, and $\overline{\bm{s}}$ the estimated factor vector.
\begin{itemize}
	\item $n,m$ - the dimensions of $\bm{A}$.
	\item $p$ - the number of features in the feature matrix.
	\item $k$ - the true number of features.
	\item $\mu$ - The fraction of missing entries in $\bm{A}$.
	\item $T$ - the total time taken for the algorithm.
	\item MAPE - the Mean Absolute Percentage Error (MAPE) for the retrieved matrix $\hat{\bm{A}}$:
	\[\text{MAPE}=\frac{1}{|\ma{S} |}\sum_{(i,j) \in  \ma{S}} \frac{|\hat{A}_{ij}-A_{ij}|}{|A_{ij}|}, \]
	where $\ma{S}=\Omega^c$ is the  set of missing  data in $\bm{A}$. 
\end{itemize}
Since the concept of Interpretable Matrix Completion is new, there is a lack of directly comparable algorithms in the literature. Thus, in lieu, we compare OptComplete to state-of-the-art solvers for Inductive Matrix Completion and general matrix completion, which are:
\begin{itemize}
	\item IMC by \cite{natarajan2014inductive} - This algorithm is a well-accepted benchmark for testing Inductive Matrix Completion algorithms. 
	\item SoftImpute-ALS (SIALS) by \cite{softimputeals} - This is widely recognized as a state-of-the-art matrix completion method without feature information. It has among the best scaling behavior across all classes of matrix completion algorithms as it utilizes fast alternating least squares to achieve  scalability. 
\end{itemize}
We use the best existing implementations of IMC (Matlab 2018b) and SIALS (R 3.4.4, package softImpute) with parallelization on the same server. 

We further compare our algorithm to CutPlanes, the original cutting plane algorithm developed in Section \ref{sec:MatrixCompletion}. It is known that for general mixed-integer convex problems, the cutting plane algorithm has the best overall performance (see e.g. \cite{lubin2016extended} for details), and thus CutPlanes represent a good baseline of comparison for OptComplete.

We randomly selected $20\%$ of those elements masked to serve as a validation set. The regularization  parameter $\gamma$ of OptComplete, the rank parameter of IMC and the penalization parameter $\lambda$ of IMC and SIALS are selected using the validation set. The results are separated into sections below. The first five sections modify one single variable out of $n,m,p,k,\mu$ to investigate OptComplete's scalability, where the leftmost column indicates the variable modified. The last section compares the  four algorithms scalability for a variety of parameters that reflect more realistic scenarios.
\vfill
\pagebreak
\begin{table}[H]
	\centering
	\resizebox{\columnwidth}{!}{%
		\begin{tabular}{|c|l|l|l|l|l||l|c|l|c|l|c|l|c|}
			\hline \multirow{2}{*}{} & \multirow{2}{*}{$\bm{n}$}&\multirow{2}{*}{$\bm{m}$}&\multirow{2}{*}{$\bm{p}$}&\multirow{2}{*}{$\bm{k}$} & \multirow{2}{*}{$\bm{\mu\%}$} & \multicolumn{2}{c|}{\textbf{OptComplete}} & \multicolumn{2}{c|}{\textbf{CutPlanes}} & \multicolumn{2}{c|}{\textbf{IMC}}& \multicolumn{2}{c|}{\textbf{SIALS}}\\\cline{7-14}
			&  &  & &  & & $\bm{T}$ & \textbf{MAPE} & $\bm{T}$ & \textbf{MAPE} & $\bm{T}$ & \textbf{MAPE}& $\bm{T}$ & \textbf{MAPE} \\\hline
			\multirow{4}{*}{$\mu$} &	100 & 100 & 15 & 5 & $20\%$ & 1.7s & $0.1\%$& 6.0s & $0.1\%$ & 0.03s & $0.01\%$& 0.02s & $0.3\%$ \\\cline{2-14}
			&100 & 100 & 15 & 5 & $50\%$ & 0.9s &$0.02\%$  & 4.5s &$0.02\%$ &0.07s & $0.5\%$ &0.03s & $0.9\%$ \\\cline{2-14}
			&100 & 100 & 15 & 5 & $80\%$ & 0.6s &$0.03\%$ & 2.5s &$0.03\%$&0.09s & $1.3\%$&0.06s & $5.6\%$ \\\cline{2-14}
			&100 & 100 & 15 & 5 & $95\%$ & 0.2s &$0.04\%$ & 1.2s &$0.04\%$&0.12s & $12.1\%$ &0.12s & $7.4\%$\\\hline\hline
			\multirow{4}{*}{$n$}&100 & 100 & 15 & 5 & $50\%$ & 0.9s &$0.02\%$ & 4.5s &$0.02\%$ &0.07s & $0.5\%$ &0.03s & $0.9\%$ \\\cline{2-14}
			&$10^3$ & 100 & 15 & 5 & $50\%$ & 3.1s &$0.01\%$& 72.5s &$0.01\%$ &0.6s & $0.4\%$ &0.1s & $0.2\%$\\\cline{2-14}
			&$10^4$ & 100 & 15 & 5 & $50\%$ & 9.5s &$0.004\%$ & 957s &$0.004\%$&4.5s & $0.3\%$&6.5s & $0.5\%$ \\\cline{2-14}
			&$10^5$ & 100 & 15 & 5 & $50\%$ & 18.0s &$0.003\%$ & 10856s &$0.003\%$ &32.7s & $0.1\%$ &38s & $3.0\%$ \\\hline\hline
			\multirow{4}{*}{$m$}&100 & 100 & 15 & 5 & $50\%$ & 0.9s &$0.02\%$ & 4.5s &$0.02\%$&0.07s & $0.5\%$ &0.03s & $0.9\%$ \\\cline{2-14}
			&100 & $10^3$ & 15 & 5 & $50\%$ & 0.7s&$0.01\%$ & 18.6s&$0.01\%$ &0.8s & $0.3\%$&0.1s & $0.5\%$ \\\cline{2-14}
			&100 & $10^4$ & 15 & 5 & $50\%$ & 1.2s &$0.004\%$ & 68.5s &$0.004\%$&6.2s &  $0.2\%$&0.8s &  $0.3\%$\\\cline{2-14}
			&100 & $10^5$ & 15 & 5 & $50\%$ & 3.0s &$0.002\%$ & 259s &$0.002\%$&56.2s & $0.1\%$ &12.7s & $0.8\%$\\\hline\hline
			\multirow{4}{*}{$p$}&100 & 100 & 15 & 5 & $50\%$ & 0.9s &$0.02\%$& 4.5s &$0.02\%$ &0.07s & $0.5\%$ &0.03s & $0.9\%$ \\\cline{2-14}
			&100 & 100 & 50 & 5 & $50\%$ & 2.0s &$0.02\%$& 18.0s &$0.02\%$ &0.3s & $0.6\%$ &0.03s & $0.9\%$\\\cline{2-14}
			&100 & 100 & 200 & 5 & $50\%$ & 12.1s &$0.02\%$& 95.9s &$0.02\%$&1.9s & $0.8\%$ &0.03s & $0.9\%$\\\cline{2-14}
			&100 & 100 & $10^3$ & 5 & $50\%$ & 90.3s &$0.02\%$& 680s &$0.02\%$ &10.4s & $1.0\%$ &0.03s & $0.9\%$ \\\hline\hline
			\multirow{4}{*}{$k$}&100 & 100 & 50 & 5 & $50\%$ & 2.0s &$0.02\%$ & 18.0s &$0.02\%$ &0.3s & $0.5\%$ &0.03s & $0.9\%$ \\\cline{2-14}
			&100 & 100 & 50 & 10 & $50\%$ & 20.7s&$0.06\%$ & 130s & 0.06\% &0.20s & $1.2\%$&0.1s & $0.8\%$\\\cline{2-14}
			&100 & 100 & 50 & 20 & $50\%$ & 240s &$0.07\%$ & 1584s &$0.07\%$ &0.35s & $2.1\%$&0.21s & $1.0\%$\\\cline{2-14}
			&100 & 100 & 50 & 30 & $50\%$ & 980s &$0.09\%$ & 8461s &$0.09\%$ &0.5s & $3.3\%$&0.43s & $2.8\%$\\\hline\hline		
			&100 & 100 & 15 & 5 & $95\%$ & 0.2s &$0.04\%$ & 1.2s &$0.04\%$ &0.12s & $12.1\%$ &0.12s & $7.4\%$ \\\cline{2-14}
			&$10^3$ & $10^3$ & 50 & 5 & $95\%$ & 1.4s &$0.006\%$ & 3.5s &$0.006\%$ &$4.6s$ & $4.7\%$&$2.8s$ & $12.5\%$ \\\cline{2-14}
			&$10^4$ & $10^3$ & 100 & 5 & $95\%$ & 5.7s &$0.002\%$  & 35.2s &$0.002\%$ &$18s$ & $2.5\%$&$20.7s$ & $12.6\%$ \\\cline{2-14}
			&$10^5$ & $10^3$ & 200 & 10 & $95\%$ & 52s &$0.001\%$ & 1520s &$0.001\%$ & 295s & $1.7\%$& $420s$ & $4.6\%$\\\cline{2-14}
			&$10^5$ & $10^4$ & 200 & 10 & $95\%$ & 98s &$0.001\%$ &5769s  & $0.001\%$ &$1750s$ &$0.5\%$&$4042s$ & $4.1\%$ \\\cline{2-14}
			&$10^6$ & $10^4$ & 200 & 10 & $95\%$ & 480s &$0.001\%$ &$N/A$ & $N/A$ & $13750s$ & $0.3\%$ &$25094s$ & $2.5\%$  \\\cline{2-14}
			&$10^6$ & $10^5$ & 200 & 10 & $95\%$ & 680s &$0.001\%$ & $N/A$ & $N/A$ & $N/A$ & $N/A$ & $N/A$ & $N/A$  \\\cline{2-14}
			&$10^6$ & $10^6$ & 200 & 10 & $95\%$ & 1415s &$0.001\%$ & $N/A$ & $N/A$ & $N/A$ & $N/A$ & $N/A$ & $N/A$  \\\hline\hline
		\end{tabular}
	}
	\caption{Comparison of OptComplete,  IMC and SIALS 
		on synthetic data.  $N/A$ means the algorithm did not complete running in 20 hours, corresponding to 72000 seconds. }
\end{table}
Overall, we see that OptComplete achieves near-exact retrieval on all datasets evaluated, and successfully recovers the factors in the ground truth. The solutions (and its error) also matches with that of CutPlanes, the standard cutting plane algorithm. The non-zero MAPE is due to the random noise added resulting in slightly perturbed coefficients.

For the realistic and large data sizes in the last panel, we see that OptComplete not only achieves near-exact retrieval, it does so while requiring considerably less time than IMC and SIALS at the same time. For $m,n$ on the scale of $n=10^6$ and $m=10^4$, OptComplete is over 20 times faster than IMC and over 40 times faster than SIALS. At the scale of $n=10^6$ and $m=10^5$, IMC and SIALS did not finish running within 20 hours, while OptComplete completed in just under 12 minutes. We also see that OptComplete achieves very significant speedups compared to the standard cutting plane algorithm - up to $60$x at the scale of $n=10^5$ and $m=10^4$. 

We analyze the scaling of OptComplete as a function of: 
\begin{enumerate}
	\item $\mu$ - The algorithm is able to retrieve the exact factors used even with $95\%$ of missing data. Furthermore, the running time decreased with increasing missing entries, consistent with the fact that is computational complexity scales with $|\Omega|$. 
	\item $n$ - The algorithm has good scalability in $n$, reflecting its $O(\sqrt{n}\log(n))$ type complexity. This allows the algorithm to support matrices with $n$ in the $10^6$ range. Its scaling behavior is superior to both IMC and SIALS. 
	\item $m$ - The algorithm scales exceptionally well in $m$. We observe that empirically the algorithm runtime seems to grow much slower than the theoretical $O(\sqrt{m}\log(m))$ dependence. A closer examination reveals that as $m$ increases, the number of cutting planes generated by Gurobi is decreasing. Qualitatively, this can be explained by a larger $m$ giving the algorithm more signal to find which $k$ features are the correct ones out of the $p$ ones. We note that such behavior is also exhibited by CutPlanes, as it roughly scales as $O(\sqrt{m})$ rather than the $O(m)$ as expected. We see that IMC and SIALS scales as $O(m)$.
	\item $p$ - The algorithm scales relatively well in $p$, which reflects the performance of the Gurobi solver. We empirically observe that Gurobi is generating roughly $O(1)-O(p)$ cutting planes. Thus, as each cutting plane is $O(p)$, we expect $O(p)-O(p^2)$ dependence, as is observed here. We note that OptComplete achieves similar scaling behavior as IMC in $p$. Note here the SIALS algorithm does not utilize feature information and thus a change in $p$ does not affect the algorithm's run speed.
	\item $k$ - The algorithm does not scale very well in $k$. We empirically observe that Gurobi solver is roughly generating $O(k)$ cutting planes and  each cutting plane has cubic dependence on $k$. It appears that SIALS and IMC almost have a linear scaling behavior. However, in most applications, such as recommendation systems or low-rank retrieval, $k$ is usually kept very low ($k\leq 30$), so this is not a particular concern. 
\end{enumerate}
\section{Real-World Experiments}
\label{sec:mc_compute_real}
In this section, we report on the performance of OptComplete on the Netflix Prize dataset \citep{bennett2007netflix}.
This dataset was released in a competition to predict ratings of customers on unseen movies, given over 10 million ratings scattered across $500,000$ people and $16,000$ movies. Thus, when presented in a matrix $\bm{A}$ where $A_{ij}$ represents the rating of individual $i$ on movie $j$, the goal is to complete the matrix $\bm{A}$ under a low-rank assumption. 

The feature matrix $\bm{B}$ of  OptComplete  is constructed using data from the TMDB Database, and covers 59 features that measure geography, popularity, top actors/actresses, box office, runtime, genre and more. The full list of 59 features is contained in  Appendix \ref{app:featurelist}.

For this experiment, we  included movies where all  59 features are available, and people who had at least   $5$ ratings present. This gives a matrix of 
$471,268$ people and $14,538$ movies. The slight reduction of size from the original data is due to the lack of features for about $2,000$ niche movies.
To observe the scalability of OptComplete, we created five data sets:
\begin{enumerate}
	\item Base - $\bm{A}_1$ has  dimensions   $3,923\times 103$.
	\item Small - $\bm{A}_2$ has  dimensions   $18,227\times 323$.
	\item Medium - $\bm{A}_3$ has  dimensions   $96,601\times 788$.
	\item Large - $\bm{A}_4$ has  dimensions   $471,268 \times 1760$.
	\item Full - $\bm{A}$ has  dimensions   $471,268\times 14,538$. 
\end{enumerate}
These sizes are constructed such that the total number of elements in $\bm{A}$ in the successive sizes are approximately different by approximately an order of magnitude.

For each individual matrix, we uniformly randomly withhold $20\%$ of the ratings as a test set $\ma{S}$, and use the remaining $80\%$ of ratings to impute a complete matrix $\hat{\bm{A}}$ - we perform cross-validation  on the appropriate hyperparameters. Then, we report MAPE.

For comparison, we  again use IMC and SIALS. We set the maximum rank of SIALS to be $k$ - the rank optimized for in OptComplete. The results are listed below:
\begin{table}[H]
	\centering
	\resizebox{\columnwidth}{!}{%
		\begin{tabular}{|l|l|l|l|l||l|l|l|l|l|l|}
			\hline  \multirow{2}{*}{$\bm{n}$}&\multirow{2}{*}{$\bm{m}$}&\multirow{2}{*}{$\bm{p}$}&\multirow{2}{*}{$\bm{k}$} & \multirow{2}{*}{$\bm{\mu\%}$} & \multicolumn{2}{c|}{\textbf{OptComplete}} & \multicolumn{2}{c|}{\textbf{IMC}}& \multicolumn{2}{c|}{\textbf{SIALS}}\\\cline{6-11}
			&  & &  & & $\bm{T}$  & \textbf{MAPE} & $\bm{T}$  & \textbf{MAPE} & $\bm{T}$  & \textbf{MAPE}\\\hline
			3,923 & 103 & 59 & 5 & $92.6\%$ & 6.0s & $29.4\%$ &0.6s & $34.2\%$& 0.3s & $31.2\%$ \\\hline
			18,227 & 323 & 59 & 5 & $94.8\%$ & 12.2s & $21.8\%$ & 5.2s & $29.1\%$ & 4.1s & $24.1\%$\\\hline
			96,601 & 788 & 59 & 5 & $94.2\%$ & 25.5s & $20.9\%$ & 38.1s & $28.7\%$ & 30.4s & $21.3\%$\\\hline
			471,268 & 1,760 & 59 & 5 & $93.6\%$ & 102s & $18.8\%$ & 460s & $24.6\%$& 430s & $19.8\%$ \\\hline
			471,268 & 14,538 & 59 & 5 & $94.1\%$ & 170s & $15.7\%$ & 3921s & $21.5\%$ &5300s & $16.7\%$\\\hline
		\end{tabular}
	}
	\caption{Comparison of methods on Netflix data for $k=5$.}
\end{table}

\begin{table}[H]
	\centering
	\resizebox{\columnwidth}{!}{%
		\begin{tabular}{|l|l|l|l|l||l|l|l|l|l|l|}
			\hline  \multirow{2}{*}{$\bm{n}$}&\multirow{2}{*}{$\bm{m}$}&\multirow{2}{*}{$\bm{p}$}&\multirow{2}{*}{$\bm{k}$} & \multirow{2}{*}{$\bm{\mu\%}$} & \multicolumn{2}{c|}{\textbf{OptComplete}} & \multicolumn{2}{c|}{\textbf{IMC}}& \multicolumn{2}{c|}{\textbf{SIALS}}\\\cline{6-11}
			&  & &  & & $\bm{T}$  & \textbf{MAPE} & $\bm{T}$  & \textbf{MAPE} & $\bm{T}$  & \textbf{MAPE}\\\hline
			3,923 & 103 & 59 & 10 & $92.6\%$ & 11.0s & $30.4\%$ &1.4s & $36.7\%$& 0.8s & $35.8\%$ \\\hline
			18,227 & 323 & 59 & 10 & $94.8\%$ & 20.3s & $24.0\%$ & 12.5s & $32.5\%$ & 7.0s & $28.9\%$\\\hline
			96,601 & 788 & 59 & 10 & $94.2\%$ & 45.9s & $22.3\%$ & 84.2s & $29.6\%$ & 50.7s & $22.8\%$\\\hline
			471,268 & 1,760 & 59 & 10 & $93.6\%$ & 260s & $20.7\%$ & 1022s & $24.8\%$& 870s & $20.7\%$ \\\hline
			471,268 & 14,538 & 59 &10  & $94.1\%$ & 380s & $19.6\%$ & 8704s & $23.1\%$ &10240s & $20.0\%$\\\hline
		\end{tabular}
	}
	\caption{Comparison of methods on Netflix data for $k=10$.}
\end{table}
We can see that OptComplete  outperforms both IMC and SIALS in accuracy across the datasets under different $k$; furthermore in the two largest datasets OptComplete ran 10x to 20x faster than IMC and SIALS. Here we see that an increase from $k=5$ to $k=10$ actually decreased out-of-sample performance as additional factors are actually not very helpful in predictive customer tastes. The decline for OptComplete and IMC were especially higher due to the fact that the possible factors are fixed and thus an increase in the number of factors caused some non-predictive factors to be included.

For the $k=5$ case, OptComplete identified the following as the top factors that influences an individual's rating:
\begin{itemize}
	\item IMDB Rating
	\item Genre: Drama
	\item Released within last 10 years
	\item Number of Top 100 Actors
	\item Produced in US
\end{itemize}
These factors provide an intuitive explanation of the individual ratings of each customer in terms of a small number of factors, while exceeding the high predictive accuracy of SIALS.

\section{Conclusions}
\label{sec_mc_conclusion}
We have presented OptComplete, a scalable algorithm  to retrieve a low-rank matrix in the presence of side information.
Compared with state of the art algorithms for matrix completion, OptComplete exceeds current benchmarks on both scalability and accuracy
and provides insight on the factors that affect the ratings. 

\clearpage

\begin{appendices}

\setcounter{table}{0}
\renewcommand{\thetable}{A\arabic{table}}
\setcounter{figure}{0}
\renewcommand{\thefigure}{A\arabic{figure}}
\setcounter{equation}{0}
\renewcommand{\theequation}{A\arabic{equation}}

\section{Proof of Theorem \ref{theo:inner_theo}}
\label{app:inner_theo_proof}
We first note that since $S_k^p$ is a finite set, we only need to prove the result for a particular $\bm{s} \in S_k^p$, and it would apply for all $\bm{s}$. Therefore, we would assume $\bm{s}$ is fixed below. For simplicity, we would only demonstrate the proof for $\alpha_i(\bm{s})=\overline{\bm{a}}_i\bm{\gamma}_i(s)$, as the one for $(\bm{B}\bm{W}_i\bm{\gamma}_i(s))^2$ follows in the same exact fashion. Furthermore, since we are only focusing on one particular $i \in \{1,\ldots,n\}$, we would drop all $i$ subscripts below for ease of notation. The quantities of interest are therefore:
\begin{align*}
    \alpha(\bm{s})=\ &\frac{1}{m}\left(\overline{\bm{a}}\overline{\bm{a}}^T
-\overline{\bm{a}}\bm{V}\left(\frac{\bm{I}_k}{\gamma}+\bm{V}^T\bm{W}\bm{V}\right)^{-1}\bm{V}^T\overline{\bm{a}}^T\right)\\
    \alpha^F(\bm{s})=\ &\frac{1}{f}\left(\overline{\bm{a}}_F\overline{\bm{a}}_F^T
-\overline{\bm{a}}_F\bm{V}_F\left(\frac{\bm{I}_k}{\gamma}+\bm{V}_F^T\bm{W}\bm{V}_F\right)^{-1}\bm{V}_F^T\overline{\bm{a}}_F^T\right).
\end{align*}

First, let $\overline{\bm{V}}=\bm{W}\bm{V}$, and let us consider a reduced QR factorization of $\overline{\bm{V}}=\bm{Q}\bm{R}$  where $\bm{Q} \in \md{R}^{m \times k}$ has orthogonal columns such that $\bm{Q}^T\bm{Q}= m \cdot \bm{I}_k$, and $\bm{R}_i \in \md{R}^{k \times k}$. Note such definition implies $\|\bm{R}\|=O(1)$. Then, we would rewrite the terms as follows:
\begin{align*}
    \alpha(\bm{s})=\ & \frac{\overline{\bm{a}}\overline{\bm{a}}^T}{m}-\frac{\overline{\bm{a}}\bm{V}}{m} \left(\frac{\bm{I}_k}{m\gamma}+\frac{\bm{R}^T\bm{Q}^T\bm{Q}\bm{R}}{m}\right)^{-1}\frac{\bm{V}^T\overline{\bm{a}}^T}{m}, \vspace{3pt} \\
    \alpha^F(\bm{s})=\ & \frac{\overline{\bm{a}}_F\overline{\bm{a}}^T_F}{f}-\frac{\overline{\bm{a}}_F\bm{V}_F}{f} \left(\frac{\bm{I}_k}{f\gamma}+\frac{\bm{R}^T\bm{Q}_F^T\bm{Q}_F\bm{R}}{f}\right)^{-1}\frac{\bm{V}_F^T\overline{\bm{a}}_F^T}{f} .
\end{align*}
We note that
\begin{alignat*}{3}
    \overline{\bm{a}}\overline{\bm{a}}^T&=\sum_{i=1}^m \overline{a}_i^2,\qquad  &\qquad  \overline{\bm{a}}\bm{V}&=\sum_{i=1}^m \overline{a}_i\bm{v}_i, \vspace{3pt} \\
    \overline{\bm{a}}_F\overline{\bm{a}}_F^T&=\sum_{i\in F} \overline{a}_i^2,\qquad  &\qquad  \overline{\bm{a}}_F\bm{V}_F&=\sum_{i\in F} \overline{a}_i\bm{v}_i.
\end{alignat*}
Therefore, if we treat $\overline{a}_1^2,\ldots \overline{a}_m^2$ as a finite population, then $\overline{\bm{a}}_F\overline{\bm{a}}_F^T$ is a random sample of $f$ points drawn without replacement from that set, and similarly for $\overline{\bm{a}}_F\bm{V}_F$. Therefore, we can then utilize Hoeffding's inequality to bound the deviation of these terms, as reproduced below:
\begin{prop}[Hoeffding's Inequality]
\label{prop:hoeffding}
Let $\mathcal{X}=(x_1,\ldots,x_n)$ be a finite population of $N$ points and $X_1,\ldots, X_n$ be a random sample drawn without replacement from $\bm{X}$. Let
\begin{equation*}
    a = \min_{1\leq i \leq n} x_i \qquad \text{and} \qquad b = \max_{1\leq i \leq n} x_i.
\end{equation*}
Then, for all $\epsilon>0$, we have
\begin{equation}
    \md{P}\left(\Big|\frac{\sum_{i=1}^n X_i}{n} -\mu \Big|\geq \epsilon\right) \leq 2\exp\left(-\frac{2n\epsilon^2}{(b-a)^2}\right).
\end{equation}
\end{prop}
For a proof, see for example \cite{boucheron2013concentration}. Then, applying Proposition \ref{prop:hoeffding} to $\overline{\bm{a}}_f\overline{\bm{a}}_f^T$, $\overline{\bm{a}}_f\bm{V}_f$, and inverting the inequality, we have
\begin{align}
    \md{P}\left(\left|\frac{\overline{\bm{a}}_F\overline{\bm{a}}^T_F}{f}-\frac{\overline{\bm{a}}\overline{\bm{a}}^T}{m}\right|\leq \sqrt{\frac{A\log(\frac{1}{\epsilon})}{f}}\right)&\geq 1- \epsilon,  \label{eq:matrix_dev_bound1}\vspace{3pt} \\
    \md{P}\left(\left\|\frac{\overline{\bm{a}}_F\bm{V}_F}{f}-\frac{\overline{\bm{a}}\bm{V}}{m}\right\|\leq \sqrt{\frac{Bk\log(\frac{k}{\epsilon})}{f}}\right)&\geq 1- \epsilon, \label{eq:matrix_dev_bound2}   
\end{align}
where $A,B$ are constants independent of $k,f,m,\epsilon$. 

Now we would show that $\frac{\bm{R}^T\bm{Q}^T\bm{Q}\bm{R}}{m}$ is close to $\frac{\bm{R}^T\bm{Q}_F^T\bm{Q}_F\bm{R}}{f}$:
\begin{lem}
\label{lem:matrix_ineq}
\begin{equation}
    \md{P}\left(\left\|\frac{\bm{R}^T\bm{Q}^T\bm{Q}\bm{R}}{m}-\frac{\bm{R}^T\bm{Q}_F^T\bm{Q}_F\bm{R}}{f}\right\|\leq \sqrt{\frac{Ck\log(\frac{k}{\epsilon})}{f}}\right)\geq 1- \epsilon. \label{eq:matrix_dev_bound3}  
\end{equation}
\end{lem}
To prove this, we would first introduce a matrix analog of the well-known Chernoff bound, the proof of which can be found in \cite{tropp2012user}:
\begin{lem}
    \label{lem:matrix_chernoff}
    Let $\mathcal{X} \in \md{R}^{k \times k}$ be a finite set of positive-semidefinite matrices, and suppose that
    \[\max_{\bm{X} \in \mathcal{X}} \lambda_{\max}(\bm{X})\leq D,\]
    where $\lambda_{\min}/\lambda_{\max}$is the minimum/maximum eigenvalue function. Sample $\{\bm{X}_1,\ldots, \bm{X}_\ell\}$ uniformly at random without replacement. Compute:
    \[\mu_{\min}:=\ell \cdot \lambda_{\min}(\md{E} \bm{X}_1)\qquad \mu_{\max}:=\ell \cdot \lambda_{\max}(\md{E} \bm{X}_1).\]
    Then,
    \begin{align*}
        \md{P}\left\{\lambda_{\min} \left(\sum_j \bm{X}_j\right)\leq (1-\delta) \mu_{\min} \right\}&\leq k \cdot \exp\left(\frac{-\delta^2\mu_{\min}}{4D}\right), \quad \text{for } \delta \in [0,1), \vspace{3pt}\\
        \md{P}\left\{\lambda_{\max} \left(\sum_j \bm{X}_j\right)\leq (1+\delta) \mu_{\max} \right\}&\leq k \cdot \exp\left(\frac{-\delta^2\mu_{\max}}{4D}\right), \quad \text{for } \delta \geq 0 .
    \end{align*}
\end{lem}
Using this lemma, we would proceed with the proof of Lemma \ref{lem:matrix_ineq}.
\proof{Proof of Lemma \ref{lem:matrix_ineq}:}
First, we note that
\begin{align*}
 \bm{Q}^T\bm{Q}&=\sum_{i=1}^m \bm{q}_i^T\bm{q}_i, \vspace{3pt} \\
 \bm{Q}^T_F\bm{Q}_F&=\sum_{i\in F} \bm{q}_i^T\bm{q}_i,
\end{align*}
where $\bm{q}_i^T\bm{q}_i \in \md{R}^{k \times k}$ rank-one positive semi-definite matrices. Therefore, we can take $\bm{Q}^T_F\bm{Q}_F$ as a random sample of size $f$ from the set $\mathcal{X}=\{\bm{q}_i^T\bm{q}_i\}_{i=1,\ldots,m}$, which satisfies the conditions in Lemma \ref{lem:matrix_chernoff} with $D=O(k)$. Furthermore, with $\mathcal{X}$, we observe that we have $\md{E} \bm{X}_1=\frac{\bm{Q}^T\bm{Q}}{m}=\bm{I}_k$, so we have
\[\lambda_{\min}(\md{E} \bm{X}_1)=\lambda_{\max}(\md{E} \bm{X}_1)=1.\]
Therefore, we apply Lemma \ref{lem:matrix_chernoff} to $\bm{Q}_F^T\bm{Q}_F$ and obtain 
\begin{align*}
    \md{P}\left\{\lambda_{\min} \left(\bm{Q}_F^T\bm{Q}_F\right)\leq (1-\delta)f\right\}&\leq k \cdot \exp\left(\frac{-\delta^2f}{kD'}\right), \vspace{3pt} \\
        \md{P}\left\{\lambda_{\max} \left(\bm{Q}_F^T\bm{Q}_F\right)\geq (1+\delta)f\right\}&\leq k \cdot \exp\left(\frac{-\delta^2f}{kD'}\right),
\end{align*}
where we set $D=\frac{kD'}{4}$ with $D'=O(1)$. Some rearrangement gives:
\begin{equation}
\label{eq:matrix_eig_ineq}
    \md{P}\left\{\lambda_{\min} \left(\frac{\bm{Q}_F^T\bm{Q}_F}{f}\right)\geq 1-\sqrt{\frac{kD'\log\left(\frac{2k}{\epsilon}\right)}{f}} \;\;  \text{and} \;\; \lambda_{\max} \left(\frac{\bm{Q}_F^T\bm{Q}_F}{f}\right)\leq 1+\sqrt{\frac{kD'\log\left(\frac{2k}{\epsilon}\right)}{f}}\right\}\geq 1 - \epsilon.
\end{equation}
Now since $\frac{\bm{Q}^T\bm{Q}}{m}=\bm{I}_k$, we have
\begin{equation}
\label{eq:matrix_eig_ineq2}
    \lambda_{\min}\left(\frac{\bm{Q}^T\bm{Q}}{m}\right)=\lambda_{\max}\left(\frac{\bm{Q}^T\bm{Q}}{m}\right)=1
\end{equation}
Combining equation (\ref{eq:matrix_eig_ineq2})
 and (\ref{eq:matrix_eig_ineq}) gives
\begin{equation}
    \md{P}\left\{\left\|\frac{\bm{Q}_F^T\bm{Q}_F}{f}-\frac{\bm{Q}^T\bm{Q}}{m}\right\|\leq \sqrt{\frac{kD'\log\left(\frac{2k}{\epsilon}\right)}{f}}\right\}\geq 1 - \epsilon.
\end{equation}
Then, we have 
\begin{equation}
    \md{P}\left\{\left\|\frac{\bm{R}^T\bm{Q}_F^T\bm{Q}_F\bm{R}}{f}-\frac{\bm{R}^T\bm{Q}^T\bm{Q}\bm{R}}{m}\right\|\leq \|\bm{R}\|^2\sqrt{\frac{kD'\log\left(\frac{2k}{\epsilon}\right)}{f}}\right\}\geq 1 - \epsilon.
\end{equation}
Taking $C=D'\|\bm{R}\|^4\log(2)$ gives the required result. (Note $\|\bm{R}\|=O(1)$ as we setup the QR decomposition to have $\bm{Q}^T\bm{Q}=O(m)$).
\QED
\endproof
With equations (\ref{eq:matrix_dev_bound1}), (\ref{eq:matrix_dev_bound2}), and (\ref{eq:matrix_dev_bound3}), we are now ready to bound  $\alpha(\bm{s})$ and $\alpha^F(\bm{s})$. We   first introduce another lemma from matrix perturbation theory (for proof, see e.g. \cite{stewart1990matrix}).
\begin{lem}
    \label{lem:matrix_inv}
    Let $\bm{A},\bm{B}$ be invertible matrices and let $\bm{B}=\bm{A}+\bm{\Delta}$. Then, we have
    \begin{equation}
        \|\bm{A}^{-1}-\bm{B}^{-1}\|\leq \|\bm{A}^{-1}\|\|\bm{B}^{-1}\|\|\bm{\Delta}\|.
    \end{equation}
\end{lem}
Then, we have
\begin{align*}
    \|\alpha(\bm{s})-\alpha^F(\bm{s})\|&\leq \left|\frac{\overline{\bm{a}}\overline{\bm{a}}^T}{m}-\frac{\overline{\bm{a}}_F\overline{\bm{a}}^T_F}{f}\right|\\&+\left|\frac{\overline{\bm{a}}\bm{V}}{m} \left(\frac{\bm{I}_k}{m\gamma}+\frac{\bm{R}^T\bm{Q}^T\bm{Q}\bm{R}}{m}\right)^{-1}\frac{\bm{V}^T\overline{\bm{a}}^T}{m}-\frac{\overline{\bm{a}}_F\bm{V}_F}{f} \left(\frac{\bm{I}_k}{f\gamma}+\frac{\bm{R}^T\bm{Q}_F^T\bm{Q}_F\bm{R}}{f}\right)^{-1}\frac{\bm{V}_F^T\overline{\bm{a}}_F^T}{f}\right|.&\intertext{Using (\ref{eq:matrix_dev_bound1}) and triangle inequality, we have}&\leq \sqrt{\frac{A\log(\frac{1}{\epsilon})}{f}} + \left|\left(\frac{\overline{\bm{a}}\bm{V}}{m}-\frac{\overline{\bm{a}}_F\bm{V}_F}{f}\right)\left(\frac{\bm{I}_k}{m\gamma}+\frac{\bm{R}^T\bm{Q}^T\bm{Q}\bm{R}}{m}\right)^{-1}\frac{\bm{V}^T\overline{\bm{a}}^T}{m}\right|\\&+\left|\frac{\overline{\bm{a}}_F\bm{V}_F}{f}\left(\left(\frac{\bm{I}_k}{m\gamma}+\frac{\bm{R}^T\bm{Q}^T\bm{Q}\bm{R}}{m}\right)^{-1}-\left(\frac{\bm{I}_k}{f\gamma}+\frac{\bm{R}^T\bm{Q}_F^T\bm{Q}_F\bm{R}}{f}\right)^{-1}\right)\frac{\bm{V}^T\overline{\bm{a}}^T}{m}\right|\\&+\left|\frac{\overline{\bm{a}}_F\bm{V}_F}{f}\left(\frac{\bm{I}_k}{f\gamma}+\frac{\bm{R}^T\bm{Q}_F^T\bm{Q}_F\bm{R}}{f}\right)^{-1}\left(\frac{\bm{V}^T\overline{\bm{a}}^T}{m}-\frac{\bm{V}_F^T\overline{\bm{a}}_F^T}{f}\right)\right|.&\intertext{Using (\ref{eq:matrix_dev_bound2}) and Lemma \ref{lem:matrix_inv}, we have}&\leq \sqrt{\frac{A\log(\frac{1}{\epsilon})}{f}}+ \sqrt{\frac{B'k\log(\frac{k}{\epsilon})}{f}}\\&+D'\left|\frac{\overline{\bm{a}}_F\bm{V}_F}{f}\left(\frac{\bm{I}_k}{m\gamma}-\frac{\bm{I}_k}{f\gamma}+\frac{\bm{R}^T\bm{Q}^T\bm{Q}\bm{R}}{m}-\frac{\bm{R}^T\bm{Q}_F^T\bm{Q}_F\bm{R}}{f}\right)\frac{\bm{V}^T\overline{\bm{a}}^T}{m}\right|\\&+\sqrt{\frac{B''k\log(\frac{k}{\epsilon})}{f}},\intertext{for some $O(1)$ constant $D'$. Note that $\|\frac{\bm{I}_k}{m\gamma}-\frac{\bm{I}_k}{f\gamma}\|\leq \frac{M}{f}$, so using (\ref{eq:matrix_dev_bound3}), we obtain 
    \begin{equation}
            \md{P}\left(\left\|\frac{\bm{I}_k}{m\gamma}-\frac{\bm{I}_k}{f\gamma}+\frac{\bm{R}^T\bm{Q}^T\bm{Q}\bm{R}}{m}-\frac{\bm{R}^T\bm{Q}_F^T\bm{Q}_F\bm{R}}{f}\right\|\leq \sqrt{\frac{C'k\log(\frac{k}{\epsilon})}{f}}\right)\geq 1- \epsilon , \label{eq:matrix_dev_bound4}
    \end{equation}
    for some new $O(1)$ constant $C'$. Then, using (\ref{eq:matrix_dev_bound4}), we have
    }&\leq \sqrt{\frac{A\log(\frac{1}{\epsilon})}{f}}+ \sqrt{\frac{B'k\log(\frac{k}{\epsilon})}{f}}+\sqrt{\frac{C''k\log(\frac{k}{\epsilon})}{f}}+\sqrt{\frac{B''k\log(\frac{k}{\epsilon})}{f}}.&\intertext{Since $k\geq 1$, we have}&\leq \sqrt{\frac{Mk\log(\frac{k}{\epsilon})}{f}},
\end{align*}
For some $O(1)$ constant $M$ with probability at least $1-4\epsilon$. Therefore, taking $\epsilon'=\frac{\epsilon}{4}$ and absorbing the extra constant into $M$ gives the requires result.
\clearpage
\section{Proof of Theorem \ref{theo:outer_theo}}
\label{app:outer_theo_proof}

We   first prove a Hoeffding-type bound as follows
\begin{prop}
\label{prop:hoeffding_general}
Let $X_1,\ldots,X_n$ be independent (but not necessarily identically distributed) random variables which satisfy 
\begin{equation*}
    \md{P}\left(|X_i-a_i|\geq t\right)\leq \exp\left(\frac{-t^2}{\sigma_i^2}\right).
\end{equation*}
Then, we have 
\begin{equation}
    \md{P}\left(\left|\frac{\sum_{i=1}^n X_i-a_i}{n}\right|\geq t\right)\leq 4\exp\left(-\frac{n^2t^2}{4\sum_{i=1}^n \sigma_i^2}\right).
\end{equation}
\end{prop}
We note that such statement differs from Hoeffding's inequality as we do not require $\md{E}[X_i]=a_i$. The proof is as follows.
\proof{Proof of Proposition \ref{prop:hoeffding_general}:}
We first introduce a lemma known as Chernoff's bounding method (proven in \cite{chernoff1952measure}):
\begin{lem}
\label{lem:chernoff_bound}
Let $Z$ be a random variable on $\md{R}$. Then for all $t>0$, we have
\begin{equation*}
    \md{P}(Z\geq t) \leq \inf_{s >0} [e^{-st}M_Z(s)],
\end{equation*}
where $M_Z(s) $ is the moment generating function of $Z$.
\end{lem}
Let  $Y_i=X_i-a_i$. Then, we have
\begin{equation}
\label{eq:sum_bound}
    \md{P}\left(\sum_{i=1}^n Y_i\geq t\right)\leq \min_{s >0} e^{-st} \prod_{i=1}^n \md{E}\left[e^{sY_i}\right] .
\end{equation}
We aim to bound $\md{E}\left[e^{sY_i}\right]$, subject to $\md{P}(|Y_i|\geq t)\leq \exp\left(\frac{-t^2}{\sigma_i^2}\right)$. Since $e^{st}$ is an increasing function of $t$, $\md{E}\left[e^{sY_i}\right]$ is maximized when $\md{P}(Y_i\geq t) = \exp\left(\frac{-t^2}{\sigma_i^2}\right)$, which results in a probability distribution function of $f_{Y_i}(y)=2\frac{t}{\sigma_i^2}\exp\left(\frac{-t^2}{\sigma_i^2}\right)$ for $Y_i$. Then, we have
\begin{align*}
    \md{E}\left[e^{sY_i}\right]&\leq \int_0^\infty 2\frac{t}{\sigma_i^2}\exp\left(st-\frac{t^2}{\sigma_i^2}\right) \d t \\&=\int^\infty_0 2\frac{t}{\sigma_i^2}\exp\left(-\left(\frac{t}{\sigma_i}-\frac{s\sigma_i}{2}\right)^2+\frac{s^2\sigma_i^2}{4}\right) \d t\\&=2\exp\left(\frac{s^2\sigma_i^2}{4}\right)\int^\infty_{0} t\exp\left(-\left(t-\frac{s\sigma_i}{2}\right)^2\right) \d t\\&=\exp\left(\frac{s^2\sigma_i^2}{4}\right)\left(\exp\left(-\frac{s^2\sigma_i^2}{4}\right)+\sqrt{\pi}\frac{s\sigma_i}{2}\left(\erf\left(\frac{s\sigma_i}{2}\right)+1\right)\right)&\intertext{where $\erf$ is the standard error function. As the error function is upper bounded by 1, the last expression is less than or equal to:}&\leq 1+\sqrt{\pi}s\sigma_i \exp\left(\frac{s^2\sigma_i^2}{4}\right)
    \\&\leq 2 \exp\left(s^2\sigma_i^2\right).
\end{align*}
We then substitute this result into (\ref{eq:sum_bound}) and obtain 
\begin{align*}
    \md{P}\left(\sum_{i=1}^n Y_i\geq t\right)&\leq \min_{s >0} 2\exp\left(-st+s^2\sum_{i=1}^n \sigma_i^2\right)&\intertext{Note that this is minimized at $s=\frac{t}{2\sqrt{\sum_{i=1}^n \sigma_i^2}}$, so we have}&\leq 2\exp\left(-\frac{t^2}{4\sum_{i=1}^n \sigma_i^2}\right).
\end{align*}
Therefore,  
\[\md{P}\left(\frac{\sum_{i=1}^n Y_i}{n}\geq t\right)\leq 2\exp\left(-\frac{n^2t^2}{4\sum_{i=1}^n \sigma_i^2}\right).\]
By applying the previous derivation to $-Y_1,\ldots, -Y_n$, we obtain
\[\md{P}\left(\frac{\sum_{i=1}^n Y_i}{n}\leq -t\right)\leq 2\exp\left(-\frac{n^2t^2}{4\sum_{i=1}^n \sigma_i^2}\right). \]
Combining the two results completes the proof. 
\QED
\endproof
Using Proposition \ref{prop:hoeffding_general}, we bound the difference between $\tilde{c}^F_G(\bm{s}_t)$ and $c(\bm{s}_t)$. We have
\begin{align*}
    |\tilde{c}^F_G(\bm{s}_t)-c(\bm{s}_t)|&=\left|\frac{1}{n} \sum_{i=1}^n \alpha_i(\bm{s}) - \frac{1}{g} \sum_{i \in G} \alpha_i^F(\bm{s})\right|\\& \leq \left|\frac{1}{n} \sum_{i=1}^n \alpha_i(\bm{s}) - \frac{1}{g} \sum_{i \in G} \alpha_i(\bm{s})\right|+ \left|\frac{1}{g} \sum_{i \in G} \alpha_i(\bm{s})- \frac{1}{g} \sum_{i \in G} \alpha_i^F(\bm{s})\right| .&\intertext{The first term can be seen as the tail bound for a random sample of size $g$ chosen without replacement from the finite set $\{\alpha_i(\bm{s})\}_{i=1,\ldots,n}$. Thus, we can apply Hoeffding's theorem in Proposition \ref{prop:hoeffding}, and obtain that with probability at least $1-\epsilon$
    \begin{equation}
        \label{eq:mean_bound}
        \left|\frac{1}{n} \sum_{i=1}^n \alpha_i(\bm{s}) - \frac{1}{g} \sum_{i \in G} \alpha_i(\bm{s})\right|\leq \sqrt{\frac{M\log\left(\frac{1}{\epsilon}\right)}{g}}.
    \end{equation}
Substituting (\ref{eq:mean_bound}) into the expression above shows that with probability at least $1-\epsilon$ we have: 
    }
    & \leq \sqrt{\frac{M\log\left(\frac{1}{\epsilon}\right)}{g}} + \left|\frac{1}{g} \sum_{i \in G} \alpha_i(\bm{s})- \frac{1}{g} \sum_{i \in G} \alpha_i^F(\bm{s})\right|.&\intertext{Note that, for any fixed set $G$, for $i,j \in G$ with $i \neq j$, we have that $\alpha_i^F(\bm{s})$ and $\alpha_j^F(\bm{s})$ are independent (as we construct $F$ separately for every $i$). Furthermore, Theorem \ref{theo:inner_theo} can be inverted to read
    \[\md{P}\left(|\alpha_i(\bm{s})-\alpha_i^F(\bm{s})|\geq t\right)\leq \exp\left(-\frac{ft^2}{Mk\log(k)}\right).\]
    Therefore, $\alpha_i^F(\bm{s})$ satisfies the conditions of Proposition \ref{prop:hoeffding_general} with $X_i=\alpha_i^F(\bm{s})$ and parameters $\sigma_i^2=\frac{Mk\log(k)}{f}$, $a_i=\alpha_i(\bm{s})$. Then, applying Proposition \ref{prop:hoeffding_general} to the second term, we have} &\leq \sqrt{\frac{M\log\left(\frac{1}{\epsilon}\right)}{g}} + \sqrt{\frac{M'k\log(\frac{k}{\epsilon})}{fg}},&\intertext{with probability $1-2\epsilon$. As $k, f\geq 1$, a loose bound is therefore}& \leq \sqrt{\frac{M''k\log\left(\frac{k}{\epsilon}\right)}{g}},
\end{align*}
with probability $1-2\epsilon$. Therefore, taking $\epsilon'=2\epsilon$ and  $A=M''\log(2)$ we have
\begin{equation}
    \label{eq:cf_deviation_bound}
    \md{P}\left(|\tilde{c}^F_G(\bm{s}_t)-c(\bm{s}_t)|\leq \sqrt{\frac{Ak\log\left(\frac{k}{\epsilon'}\right)}{g}}\right) \geq 1- \epsilon' .
\end{equation}
Through a similar derivation, we have
\begin{equation}
    \label{eq:cf_deriv_deviation_bound}
    \md{P}\left(\|\nabla \tilde{c}^F_G(\bm{s}_t)-\nabla c(\bm{s}_t)\|_2\leq \sqrt{\frac{B(p+k)\log\left(\frac{k}{\epsilon'}\right)}{g}}\right) \geq 1- \epsilon' .
\end{equation}
\clearpage
\section{Proof of Proposition \ref{prop:strong_convexity_cond}}
\label{app:strong_convexity_cond_proof}
We would prove the contrapositive. Assume that $a=0$. We would calculate the second derivative of the cost function expression in \ref{eq:mc_central}. First, define $\bm{T}_{ij}=\frac{\bm{I}}{\gamma} + \bm{W}_i\bm{K}_j\bm{W}_i$. Then our cost function can be rewritten as:
\[c(\bm{s})=\frac{1}{nm} \sum_{i=1}^n \bar{\bm{a}}_i \left(\sum_{j=1}^p s_j \bm{T}_{ij}\right)^{-1}\bar{\bm{a}}_i^T\]
Then the second derivative can be easily calculated as:
\[\frac{\partial c(\bm{s})}{\partial s_k\partial s_l} = 2\sum_{i=1}^n\bar{\bm{a}}_i \left(\sum_{j=1}^p s_j \bm{T}_{ij}\right)^{-1}\bm{T}_{ik}\left(\sum_{j=1}^p s_j \bm{T}_{ij}\right)^{-1}\bm{T}_{il}\left(\sum_{j=1}^p s_j \bm{T}_{ij}\right)^{-1} \bar{\bm{a}}_i^T\]
Thus, for any vector $\bm{t} \in \md{R}^p$, we have:
\[\bm{t}^T\frac{\partial c(\bm{s})}{\partial s_k\partial s_l}\bm{t}= 2\sum_{i=1}^n\bar{\bm{a}}_i \left(\sum_{j=1}^p s_j \bm{T}_{ij}\right)^{-1}\left(\sum_{k=1}^p t_k\bm{T}_{ik}\right)\left(\sum_{j=1}^p s_j \bm{T}_{ij}\right)^{-1}\left(\sum_{l=1}^p t_l\bm{T}_{il}\right)\left(\sum_{j=1}^p s_j \bm{T}_{ij}\right)^{-1} \bar{\bm{a}}_i^T\]
Now define $\bm{v}_i=\left(\sum_{l=1}^p t_l\bm{T}_{il}\right)\left(\sum_{j=1}^p s_j \bm{T}_{ij}\right)^{-1} \bar{\bm{a}}_i^T$, then our expression becomes:
\[\bm{t}^T\frac{\partial c(\bm{s})}{\partial s_k\partial s_l}\bm{t}= 2\sum_{i=1}^n \bm{v}^T_i\left(\sum_{j=1}^p s_j \bm{T}_{ij}\right)^{-1}\bm{v}_i\]
Now since $a=0$, the Hessian cannot be positive definite for all $\bm{s}$ (as if it was, then $a>0$). Then there exist $\bm{s}^0,\bm{s}^1 \in S^p_k$ such that for $\bm{t}=\bm{s}^0-\bm{s}^1$, we have:
\[\bm{t}^T\frac{\partial c(\bm{s})}{\partial s_k\partial s_l} \bigg\rvert_{\bm{s}=\bm{s}^0}\bm{t}= 2\sum_{i=1}^n \bm{v}^T_i\left(\sum_{j=1}^p s_j^0 \bm{T}_{ij}\right)^{-1}\bm{v}_i=0\]
and $\bm{t}$ 
Note that for any $\gamma>0$, $\bm{T}_{ij}$ is positive definite, so $\sum_{j=1}^p s_j \bm{T}_{ij}$ is positive definite for any $\bm{s} \in S^p_k$ and all $i$. Therefore if $a=0$, we must have $\bm{v}_i=0$ for all $i$, which means that:
\[\left(\sum_{l=1}^p t_l\bm{T}_{il}\right)\left(\sum_{j=1}^p s_j \bm{T}_{ij}\right)^{-1} \bar{\bm{a}}_i^T=\bm{0}\]
for all $i$. Now note that we have:
\[c(\bm{s}^1)=c(\bm{s}^0)+\nabla c(\bm{s}^0)^T\bm{t}+ \bm{t}^T\frac{\partial c(\bm{s})}{\partial s_k\partial s_l} \bigg\rvert_{\bm{s}=\bm{s}^0}\bm{t} + \cdots\]
We have that:
\[\nabla c(\bm{s}^0)^T\bm{t}=\sum_{i=1}^n \left(\sum_{l=1}^p t_l\bm{T}_{il}\right)\left(\sum_{j=1}^p s_j \bm{T}_{ij}\right)^{-1} \bar{\bm{a}}_i^T=0\]
And similarly, we have $\bm{t}^T\frac{\partial c(\bm{s})}{\partial s_k\partial s_l} \bigg\rvert_{\bm{s}=\bm{s}^0}\bm{t}=0$ and all higher derivatives being 0.

Thus, we have:
\[c(\bm{s}^1)=c(\bm{s}^0)\]
Therefore, if we have $a=0$, then we must have equality on the cost function for two feasible solutions. Thus, if the cost function is not degenerate for feasible solutions, then $a>0$.
\clearpage
\section{Proof of Theorem \ref{theo:mc_proof}}
\label{app:optcomplete_theo_proof}
In this section, we provide the proof of Theorem \ref{theo:mc_proof}. We first note that OptComplete is a specific implementation of the outer approximation algorithm, which \cite{fletcher1994solving} have shown to always terminate in finite number of steps $C$. Thus, given that we assumed the problem is feasible, for OptComplete to not return an optimal solution, it would have to cut it off during the course of its execution. Let $\bm{s}^*$ be an optimal  solution for  Problem (\ref{mainproblem}). Let $\bm{s}_t$ be an optimal solution at the $t$-th iteration of OptComplete, $t\in[C]$. The cutting plane constraint for OptComplete
at  the point of an optimal solution $\bm{s}^*$
is 
$$ \eta_t\geq \tilde{c}_G^F(\bm{s}_t)+\nabla \tilde{c}_G^F(\bm{s}_t)^T(\bm{s}^*-\bm{s}_t).$$
If $c(\bm{s}^*)<\eta_t$, then $\bm{s}^*$ will be cut off, and OptComplete will not find $\bm{s}^*$. 
Applying the definition of the  convexity parameter (\ref{convexityparam}) and letting $\|\bm{s}^*-\bm{s}_t\|=\theta_t$ (noting that $\theta_t \geq 1$)
we obtain 
\begin{equation}
c(\bm{s}^*)\geq c(\bm{s}_t)+\nabla c(\bm{s}_t)^T(\bm{s}^*-\bm{s}_t)+\frac{\theta_t^2 a^2}{2} \label{exactcutineq}.
\end{equation}
Therefore, if 	
$$c(\bm{s}_t)+\nabla c(\bm{s}_t)^T(\bm{s}^*-\bm{s}_t)+\frac{\theta_t^2 a^2 }{2} \leq c(\bm{s}^*) <
\tilde{c}_r(\bm{s}_t)+\nabla \tilde{c}_r(\bm{s}_t)^T(\bm{s}^*-\bm{s}_t),$$ 
or equivalently if 
\begin{equation}
\zeta_t :=[\tilde{c}_G^F(\bm{s}_t)-c(\bm{s}_t)]+[\nabla \tilde{c}_G^F(\bm{s}_t) - \nabla c(\bm{s}_t)]^T(\bm{s}^*-\bm{s}_t) >  \frac{\theta_t^2 a^2}{2} \label{cutoffoptim},
\end{equation}	
then  OptComplete will not find $\bm{s}^*$. Therefore, for  OptComplete to succeed, $\zeta_t$ should satisfy $\zeta_t \leq \theta_t^2  a^2/2$ for all $t\in [C].$

Let $S=$ the event that   OptComplete  succeeds in finding $\bm{s}^*$. Then,
$$\md{P}(S)\geq \md{P}\left( \zeta_t \leq  \frac{\theta_t^2 a^2}{2}\; \; t\in[C]\right).$$
Since at each step of OptComplete we 
randomly sample $r$ new rows and $s$ new columns, the events $\zeta_t \leq  \frac{\theta_t^2 a^2}{2}$ are independent for different $t\in[C]$,
and hence
\begin{equation}
\label{eq:success_prob}
\md{P}(S)\geq \prod_{t=1}^C \left(1-\md{P}\left( \zeta_t > \frac{\theta_t^2 a^2}{2}\right)\right) .
\end{equation}
Therefore, we focus on calculating $\md{P}( \zeta_t \geq  \theta_t^2 a^2/2)$. Since
\begin{equation*}
\zeta_t :=[\tilde{c}_G^F(\bm{s}_t)-c(\bm{s}_t)]+[\nabla \tilde{c}_G^F(\bm{s}_t) - \nabla c(\bm{s}_t)]^T(\bm{s}^*-\bm{s}_t),
\end{equation*}
using Theorem \ref{theo:outer_theo}, we can thus provide the following bound for deviation of $\zeta_t$, for some constant $C$:
\begin{equation}
    \label{eq:zeta_deviation_bound}
    \md{P}\left(\zeta_t\leq \theta_t\sqrt{\frac{C(p+k)\log\left(\frac{k}{\epsilon}\right)}{g}}\right) \geq 1- \epsilon ,
\end{equation}
where $\theta_t=\|\bm{s}^*-\bm{s}_t\|$. Then we can invert this to calculate $\md{P}( \zeta_t \geq  \theta_t^2 a^2/2)$
\begin{equation}
\label{eq:zeta_deviation_prob}
    \md{P}\left(\zeta_t \geq  \theta_t^2 a^2/2\right) \leq k\exp\left(-\frac{a^4g}{4C(p+k)\theta_t^2} \right)\leq k\exp\left(-\frac{Da^4g}{(p+k)} \right),
\end{equation}
taking $D=\frac{1}{4C}$, and noting that $\theta_t\geq 1$. Then, we substitute (\ref{eq:zeta_deviation_prob}) into (\ref{eq:success_prob}) to obtain
\begin{align*}
    \md{P}(S) & \geq \left(1-k\exp\left(-\frac{Da^4g}{(p+k)}\right)\right)^C \\
    & \geq 1-kC\exp\left(-\frac{Da^4g}{(p+k)}\right).
\end{align*}
completing the proof.

\clearpage

\section{List of Features Used in the Netflix Problem}
\label{app:featurelist}
\begin{itemize}
	\item 24 Indicator Variables for Genres: Action, Adventure, Animation, Biography, Comedy, Crime, Documentary, Drama, Family, Fantasy, Film Noir, History, Horror, Music, Musical, Mystery, Romance, Sci-Fi, Short, Sport, Superhero, Thriller, War, Western
	\item 5 Indicator Variables for Release Date: Within last 10 years, Between 10-20 years, Between 20-30 years, Between 30-40 years, Between 40-50 Years
	\item 6 Indicator Variables for Top Actors/Actresses defined by their Influence Score at time of release: Top 100 Actors, Top 100 Actresses, Top 250 Actors, Top 250 Actresses, Top 1000 Actors, Top 1000 Actresses
	\item IMDB Rating
	\item Number of Reviews
	\item Total Production Budget
	\item Total Runtime
	\item Total Box Office Revenue
	\item Indicator Variable for whether it is US produced
	\item 11 Indicator Variables for Month of Year Released (January removed to prevent multicollinearity)
	\item Number of Original Music Score
	\item Number of Male Actors
	\item Number of Female Factors
	\item 3 Indicator Variables for Film Language: English, French, Japanese
	\item Constant 
\end{itemize}
\end{appendices}
\clearpage

\bibliographystyle{informs2014} 
\bibliography{IMCbib} 

\begin{thebibliography}{38}
\providecommand{\natexlab}[1]{#1}
\providecommand{\url}[1]{\texttt{#1}}
\providecommand{\urlprefix}{URL }

\bibitem[{Beck \protect\BIBand{} Teboulle(2009)}]{beck2009fast}
Beck A, Teboulle M (2009) A fast iterative shrinkage-thresholding algorithm for
  linear inverse problems. \emph{SIAM Journal on Imaging Sciences}
  2(1):183--202.

\bibitem[{Bennett et~al.(2007)Bennett, Lanning et~al.}]{bennett2007netflix}
Bennett J, Lanning S, et~al. (2007) The netflix prize. \emph{KDD Cup and
  Workshop 2007} (Citeseer).

\bibitem[{Bertsimas \protect\BIBand{} Copenhaver(2018)}]{bc18}
Bertsimas D, Copenhaver MS (2018) Characterization of the equivalence of
  robustification and regularization in linear, median, and matrix regression.
  \emph{European Journal of Operations Research} 270:931--942.

\bibitem[{Bertsimas \protect\BIBand{} van Parys(2020)}]{SparseReg}
Bertsimas D, van Parys B (2020) Sparse high-dimensional regression: Exact
  scalable algorithms and phase transitions. \emph{Annals of Statistics}
  48(1):300--323.

\bibitem[{Boucheron et~al.(2013)Boucheron, Lugosi, \protect\BIBand{}
  Massart}]{boucheron2013concentration}
Boucheron S, Lugosi G, Massart P (2013) \emph{Concentration inequalities: A
  nonasymptotic theory of independence} (Oxford university press).

\bibitem[{Boyd et~al.(1994)Boyd, El~Ghaoui, Feron, \protect\BIBand{}
  Balakrishnan}]{boyd1994linear}
Boyd S, El~Ghaoui L, Feron E, Balakrishnan V (1994) \emph{Linear matrix
  inequalities in system and control theory}, volume~15 (SIAM).

\bibitem[{Candes \protect\BIBand{} Plan(2010)}]{candes2010matrix}
Candes EJ, Plan Y (2010) Matrix completion with noise. \emph{Proceedings of the
  IEEE} 98(6):925--936.

\bibitem[{Cand{\`e}s \protect\BIBand{} Tao(2010)}]{candesexact}
Cand{\`e}s EJ, Tao T (2010) The power of convex relaxation: Near-optimal matrix
  completion. \emph{IEEE Transactions on Information Theory} 56(5):2053--2080.

\bibitem[{Chen et~al.(2014{\natexlab{a}})Chen, Bhojanapalli, Sanghavi,
  \protect\BIBand{} Ward}]{chen2014coherent}
Chen Y, Bhojanapalli S, Sanghavi S, Ward R (2014{\natexlab{a}}) Coherent matrix
  completion. \emph{International Conference on Machine Learning}, 674--682.

\bibitem[{Chen et~al.(2014{\natexlab{b}})Chen, Jalali, Sanghavi,
  \protect\BIBand{} Xu}]{chen2014clustering}
Chen Y, Jalali A, Sanghavi S, Xu H (2014{\natexlab{b}}) Clustering partially
  observed graphs via convex optimization. \emph{The Journal of Machine
  Learning Research} 15(1):2213--2238.

\bibitem[{Chernoff(1952)}]{chernoff1952measure}
Chernoff H (1952) A measure of asymptotic efficiency for tests of a hypothesis
  based on the sum of observations. \emph{The Annals of Mathematical
  Statistics} 23(4):493--507.

\bibitem[{Chiang et~al.(2015)Chiang, Hsieh, \protect\BIBand{}
  Dhillon}]{chiang2015matrix}
Chiang KY, Hsieh CJ, Dhillon IS (2015) Matrix completion with noisy side
  information. \emph{Advances in Neural Information Processing Systems},
  3447--3455.

\bibitem[{Chiang et~al.(2014)Chiang, Hsieh, Natarajan, Dhillon,
  \protect\BIBand{} Tewari}]{chiang2014prediction}
Chiang KY, Hsieh CJ, Natarajan N, Dhillon IS, Tewari A (2014) Prediction and
  clustering in signed networks: a local to global perspective. \emph{The
  Journal of Machine Learning Research} 15(1):1177--1213.

\bibitem[{Dhillon et~al.(2013)Dhillon, Lu, Foster, \protect\BIBand{}
  Ungar}]{dhillon2013new}
Dhillon P, Lu Y, Foster DP, Ungar L (2013) New subsampling algorithms for fast
  least squares regression. \emph{Advances in Neural Information Processing
  Systems}, 360--368.

\bibitem[{Duran \protect\BIBand{} Grossmann(1986)}]{CuttingPlane}
Duran MA, Grossmann IE (1986) An outer-approximation algorithm for a class of
  mixed-integer nonlinear programs. \emph{Mathematical Programming}
  36(3):307--339.

\bibitem[{Farhat et~al.(2013)Farhat, Shapiro, Kieser, Sultana, Jacobson,
  Victor, Warren, Streicher, Calver, Sloutsky et~al.}]{farhat2013genomic}
Farhat MR, Shapiro BJ, Kieser KJ, Sultana R, Jacobson KR, Victor TC, Warren RM,
  Streicher EM, Calver A, Sloutsky A, et~al. (2013) Genomic analysis identifies
  targets of convergent positive selection in drug-resistant mycobacterium
  tuberculosis. \emph{Nature Genetics} 45(10):1183.

\bibitem[{Fletcher \protect\BIBand{} Leyffer(1994)}]{fletcher1994solving}
Fletcher R, Leyffer S (1994) Solving mixed integer nonlinear programs by outer
  approximation. \emph{Mathematical programming} 66(1-3):327--349.

\bibitem[{Hastie et~al.(2015)Hastie, Mazumder, Lee, \protect\BIBand{}
  Zadeh}]{softimputeals}
Hastie T, Mazumder R, Lee JD, Zadeh R (2015) Matrix completion and low-rank svd
  via fast alternating least squares. \emph{The Journal of Machine Learning
  Research} 16(1):3367--3402.

\bibitem[{Jain \protect\BIBand{} Dhillon(2013)}]{jain2013provable}
Jain P, Dhillon IS (2013) Provable inductive matrix completion. \emph{arXiv
  preprint arXiv:1306.0626} .

\bibitem[{Jain et~al.(2010)Jain, Meka, \protect\BIBand{}
  Dhillon}]{jain2010guaranteed}
Jain P, Meka R, Dhillon IS (2010) Guaranteed rank minimization via singular
  value projection. \emph{Advances in Neural Information Processing Systems},
  937--945.

\bibitem[{Ji et~al.(2010)Ji, Liu, Shen, \protect\BIBand{} Xu}]{ji2010robust}
Ji H, Liu C, Shen Z, Xu Y (2010) Robust video denoising using low rank matrix
  completion. \emph{Computer Society Conference on Computer Vision and Pattern
  Recognition} (IEEE).

\bibitem[{Keshavan et~al.(2009)Keshavan, Oh, \protect\BIBand{}
  Montanari}]{keshavan2009matrix}
Keshavan RH, Oh S, Montanari A (2009) Matrix completion from a few entries.
  \emph{Information Theory, 2009. ISIT 2009. IEEE International Symposium on},
  324--328 (IEEE).

\bibitem[{Koren et~al.(2009)Koren, Bell, \protect\BIBand{}
  Volinsky}]{koren2009matrix}
Koren Y, Bell R, Volinsky C (2009) Matrix factorization techniques for
  recommender systems. \emph{Computer} 30--37.

\bibitem[{Lu et~al.(2016)Lu, Liang, Sun, \protect\BIBand{} Bi}]{lu2016sparse}
Lu J, Liang G, Sun J, Bi J (2016) A sparse interactive model for matrix
  completion with side information. \emph{Advances in Neural Information
  Processing Ssystems}, 4071--4079.

\bibitem[{Lubin et~al.(2016)Lubin, Yamangil, Bent, \protect\BIBand{}
  Vielma}]{lubin2016extended}
Lubin M, Yamangil E, Bent R, Vielma JP (2016) Extended formulations in
  mixed-integer convex programming. \emph{International Conference on Integer
  Programming and Combinatorial Optimization}, 102--113 (Springer).

\bibitem[{Mazumder et~al.(2010)Mazumder, Hastie, \protect\BIBand{}
  Tibshirani}]{SoftImpute}
Mazumder R, Hastie T, Tibshirani R (2010) Spectral regularization algorithms
  for learning large incomplete matrices. \emph{Journal of Machine Learning
  Research} 11(Aug):2287--2322.

\bibitem[{Natarajan \protect\BIBand{} Dhillon(2014)}]{natarajan2014inductive}
Natarajan N, Dhillon IS (2014) Inductive matrix completion for predicting
  gene--disease associations. \emph{Bioinformatics} 30(12):160--168.

\bibitem[{Nazarov et~al.(2018)Nazarov, Shirokikh, Burkina, Fedonin,
  \protect\BIBand{} Panov}]{nazarov2018sparse}
Nazarov I, Shirokikh B, Burkina M, Fedonin G, Panov M (2018) Sparse group
  inductive matrix completion. \emph{arXiv preprint arXiv:1804.10653} .

\bibitem[{Negahban \protect\BIBand{} Wainwright(2012)}]{negahban2012restricted}
Negahban S, Wainwright MJ (2012) Restricted strong convexity and weighted
  matrix completion: Optimal bounds with noise. \emph{Journal of Machine
  Learning Research} 13(May):1665--1697.

\bibitem[{Recht \protect\BIBand{} R{\'e}(2013)}]{recht2013parallel}
Recht B, R{\'e} C (2013) Parallel stochastic gradient algorithms for
  large-scale matrix completion. \emph{Mathematical Programming Computation}
  5(2):201--226.

\bibitem[{Shah et~al.(2017)Shah, Rao, \protect\BIBand{} Ding}]{shah2017matrix}
Shah V, Rao N, Ding W (2017) Matrix factorization with side and higher order
  information. \emph{Stat} 1050:4.

\bibitem[{Si et~al.(2016)Si, Chiang, Hsieh, Rao, \protect\BIBand{}
  Dhillon}]{si2016goal}
Si S, Chiang KY, Hsieh CJ, Rao N, Dhillon IS (2016) Goal-directed inductive
  matrix completion. \emph{Proceedings of the 22nd ACM SIGKDD International
  Conference on Knowledge Discovery and Data Mining}, 1165--1174 (ACM).

\bibitem[{Soni et~al.(2016)Soni, Chevalier, \protect\BIBand{}
  Jain}]{soni2016noisy}
Soni A, Chevalier T, Jain S (2016) Noisy inductive matrix completion under
  sparse factor models. \emph{arXiv preprint arXiv:1609.03958} .

\bibitem[{Stewart(1990)}]{stewart1990matrix}
Stewart GW (1990) Matrix perturbation theory .

\bibitem[{Tanner \protect\BIBand{} Wei(2013)}]{tanner2013normalized}
Tanner J, Wei K (2013) Normalized iterative hard thresholding for matrix
  completion. \emph{SIAM Journal on Scientific Computing} 35(5):S104--S125.

\bibitem[{Tropp(2012)}]{tropp2012user}
Tropp JA (2012) User-friendly tail bounds for sums of random matrices.
  \emph{Foundations of computational mathematics} 12(4):389--434.

\bibitem[{Woodbury(1949)}]{matrixinv}
Woodbury MA (1949) The stability of out-input matrices. \emph{Chicago, Ill} .

\bibitem[{Xu et~al.(2013)Xu, Jin, \protect\BIBand{} Zhou}]{xuspeedup}
Xu M, Jin R, Zhou ZH (2013) Speedup matrix completion with side information:
  Application to multi-label learning. \emph{Advances in Neural Information
  Processing Systems}, 2301--2309.

\end{thebibliography}

\end{document}